\documentclass[a4paper,11pt]{article}
\usepackage[english]{babel}

\usepackage[dvips]{graphicx}

\usepackage[hypertex]{hyperref}
\usepackage{latexsym}
\usepackage{amsfonts}
\usepackage{amsmath}
\usepackage{amsthm}
\usepackage{amssymb}
\usepackage{geometry} 
\usepackage{subfigure}
\usepackage{changebar}
\outerbarstrue
\usepackage{eucal}

 \theoremstyle{plain}
    \newtheorem{theorem}{Theorem}
    \newtheorem{lemma}{Lemma}
    \newtheorem{proposition}{Proposition}

\newcommand{\mathbbm}{\mathbb}
\newcommand{\ie}{ {\it i.e.~}}

\newcommand{\conf}{\mathcal{C}}
\newcommand{\pat}{\mathcal{P}} 
\newcommand{\e}[1]{\mathbf{#1}} 
\newcommand{\K}{\mathsf{K}}  
\newcommand{\Ki}{\mathsf{K}^{\scriptscriptstyle -1}} 
\newcommand{\bv}{\mathbf{b}} 
\newcommand{\wv}{\mathbf{w}} 

\newcommand{\esp}[1]{\mathbbm{E}\left[#1\right]} 
\newcommand{\espind}[2]{\mathbbm{E}_{#1}\left[#2\right]} 
\newcommand{\prob}[1]{\mathbbm{P}\left[#1\right]} 

\newcommand{\ud}{\mathrm{d}}                     
\newcommand{\xx}{\mathbf{x}}
\newcommand{\flucN}[1]{\tilde{N}_{#1}^\varepsilon}
\newcommand{\lflucN}[1]{\mathcal{N}_{#1}}

\newcommand{\tore}{\mathbbm{T}^2}
\newcommand{\cercle}{\mathbbm{S}^1}
\newcommand{\RR}{\mathbbm{R}}
\newcommand{\ZZ}{\mathbbm{Z}}
\newcommand{\CC}{\mathbbm{C}}

\newcommand{\epsi}{\varepsilon}

\newcommand{\smallterms}{\begin{array}{c} \textrm{small} \\ \textrm{terms} \end{array}}
\newcommand{\oscterms}{\begin{array}{c} \textrm{oscillating} \\ \textrm{terms} \end{array}}

\renewcommand{\Re}{\mathrm{Re}}
\renewcommand{\Im}{\mathrm{Im}}
\renewcommand{\div}{\mathrm{div}}

\DeclareMathOperator{\tr}{tr}

\DeclareMathOperator{\sgn}{sgn}

\DeclareMathOperator{\Cof}{Cof}  

\title{Pattern densities in fluid dimer models}
\author{C\'edric Boutillier\thanks{
Centrum voor Wiskunde en Informatica - P.O. Box 94079, NL-1090 GB Amsterdam\newline \texttt{cedric.boutillier@cwi.nl}}}

\begin{document}

\maketitle

\begin{abstract}
In this paper, we introduce a family of observables for the dimer model on a bi-periodic bipartite planar graph, called \emph{pattern density fields}. We study the scaling limit of these objects for liquid and gaseous Gibbs measures of the dimer model, and prove that they converge to a linear combination of a derivative of the Gaussian massless free field and an independent white noise.
\end{abstract}

\section{Introduction}
\subsection{The dimer model}
Let $G$ be a graph. 
A \emph{dimer configuration} $\conf$ of a graph $G$ is a subset of edges of $G$ such that every vertex of $G$ is incident to exactly one edge of $\conf$. The \emph{dimer model} is a system from statistical mechanics, obtained by endowing the set of all possible dimer configurations with a probability measure. It has been introduced in the 1930s~\cite{FowlerRushbrooke} to give a model for adsorption of diatomic molecules (dimers) on the surface of the crystal, represented by the graph $G$. This model of statistical mechanics is one of the rare models that can be solved exactly. For an  introduction to the dimer model, see for example~\cite{Ke:Intro}.

We will require $G$ to have at least one dimer configuration. Suppose for the moment that $G$ is finite. One can define a \emph{Boltzmann probability measure} on the dimer configurations of $G$ as follows : positive weights $w_{\e{e}}$ are assigned to the edges $\e{e}$ of $G$ and the probability of a dimer configuration $\conf$ is chosen to be proportional to the product of the weights of the edges it contains:
\begin{equation}
\prob{\conf}=\frac{1}{Z}\prod_{\e{e}\in\conf} w_{\e{e}},
\end{equation}
where the normalizing factor $Z=\sum_{\conf}\prod_{e\in\conf} w_{\e{e}}$ is called the \emph{partition function}.

Kasteleyn showed~\cite{Ka:Crys} that if $G$ is planar, the partition function $Z$ can be expressed as the Pfaffian of a weighted adjacency matrix for a well-chosen orientation of the graph. If moreover the graph is bipartite, $Z$ reduces to the determinant of a certain matrix $\K$, called the Kasteleyn operator, a cousin of the adjacency matrix of $G$, whose lines are indexed by white vertices, and columns by black vertices. In particular, it has the property that if $\e{e}=(\wv,\bv)$, $|\K(\wv,\bv)|=w[\e{e}]$. See~\cite{KOS} to see  how to define $\K$. Because of the correspondence between determinants of submatrices of a matrix and those of its inverse, the correlations are given by determinants of submatrices of $\Ki$~\cite{Ke:LocStat}.

When the graph $G$ is infinite, it may have an infinite number of dimer configurations, and it is not possible anymore to define directly a Boltzmann measure. For planar bipartite $\ZZ^2$-periodic graphs, endowed with periodic weights on edges, this notion is replaced by that of \emph{Gibbs probability measure} that has the following properties:
\begin{itemize}
\item it is ergodic under the action of $\ZZ^2$ by translation,
\item if the dimer configuration is fixed in an annular region, then the random dimer configuration inside and outside the annulus are independent, and the induced probability measure inside the annulus is the Boltzmann measure defined above.
\end{itemize}

Sheffield proved~\cite{Shef} that for given periodic weights, there exists a family of Gibbs measures leading to the same Boltzmann measure on finite regions parametrized by the two components of an external magnetic field $B$. We will only consider here the case when $B=0$, since the other measures of this family can be obtained as measures without external field for different weights.

The fundamental domain $G_1$ of $G$ is supposed to have the same number of white and black vertices. The white (resp. black) vertices are labeled, say, from 1 to $n$: $\wv_1,\dots,\wv_n$ (resp. $\bv_1,\dots,\bv_n$). If $\e{v}$ is a vertex of $G$, then $\e{v}_{\xx}=\e{v}+(x,y)$ denotes the translate of $\e{v}$ by the lattice vector $\xx=(x,y)\in \ZZ^2$. If $(z,w)$ is in the unit torus $\tore=\cercle\times\cercle$, a function $f$ is said to be $(z,w)$-periodic if
\begin{equation*}
\forall\ \e{v}\in G,\quad \forall\ (x,y)\in\ZZ^2,\quad f(\e{v}_{x,y})=z^{-y} w^x f(\e{v}).
\end{equation*}

As in the case of a planar finite region, a infinite Kasteleyn matrix $\K$ is defined. Its action on functions defined on black vertices returns a function on white vertices
\begin{equation*}
(\K f)_{\wv} = \sum_{\wv\sim\bv} \K(\wv,\bv) f_{\bv}
\end{equation*}
If $\e{e}=(\wv,\bv)$, we will often write $\K_{\e{e}}$ instead of $\K(\wv,\bv)$.

Since $\K$ is periodic, its Fourier transform $K(z,w)$, being the restriction of $\K$ to the $(z,w)$-periodic functions, is a $n\times n$ matrix. The lines of $K(z,w)$ are indexed by the white vertices of the fundamental domain $G_1$ (more precisely, by the $(z,w)$-periodic function taking value 1 at a given white vertex of the fundamental domain, and 0 elsewhere). Its columns are indexed by the black vertices of $G_1$.

The correlations between dimers under the Gibbs measure are expressed in terms of determinants whose entries are coefficients of the inverse Kasteleyn operator $\Ki$ given by the usual inverse Fourier Transform formula. Let $\bv$ and $\wv$ be respectively a black and a white vertex in the fundamental domain, and $\xx=(x,y),\xx'=(x',y')$ be lattice points. The coefficient $\Ki(\bv_{\xx},\wv_{\xx'})$ is given by
\begin{align*}
\Ki(\bv_{\xx},\wv_{\xx'})=\Ki(\bv_{\xx-\xx'},\wv)&=\iint_{\tore} z^{-(y-y')} w^{(x-x')} K^{-1}_{\bv\wv}(z,w) \frac{\ud z}{2i\pi z}\frac{\ud w}{2i\pi w}\\
&=\iint_{\tore} \frac{z^{-(y-y')} w^{(x-x')} Q_{\bv\wv}(z,w)}{P(z,w)} \frac{\ud z}{2i\pi z}\frac{\ud w}{2i\pi w}
\end{align*}
where $K^{-1}(z,w)$ is the inverse of the $n\times n$ matrix $K(z,w)$, $P(z,w)$ its determinant and $Q_{\bv,\wv}(z,w)$ is the cofactor associated to $K_{\wv,\bv}(z,w)$. 

The probability that edges $\e{e}_1=(\wv^{1},\bv^{1}), \dots \e{e}_k=(\wv^{k},\bv^{k})$ appear in the random dimer configuration is given by
\begin{equation}\label{eq:correledges}
\prob{\e{e}_1,\dots \e{e}_k}=\left(\prod_{j=1}^k \K(\wv^j,\bv^j)\right) \det_{1\leq i,j\leq k}\bigl[\Ki(\bv^{i},\wv^{j})\bigr].
\end{equation}

The way the correlation between edges decays with the distance depends on the number of zeros of $P(z,w)$ on the unit torus~\cite{KOS}. Three behaviours, or \emph{phases}, are possible for a Gibbs measure. The measure can be:
\begin{itemize}
\item \emph{solid}, when there are bi-infinite duals path on which edges appear with probability 0 or 1.
\item \emph{liquid} when $\Ki(\bv_\xx,\wv)$ decays linearly with $|\xx|$, and thus correlations decay polynomially. A more precise statement is given in lemma \ref{lem:Kiasymp1}.
\item \emph{gaseous} when $\Ki(\bv_\xx,\wv)$ and thus correlations, decay exponentially.
\end{itemize}

We will here deal with \emph{fluid} Gibbs measures, corresponding to the last two cases.

\subsection{Scaling limits of pattern densities}

A dimer configuration on a bipartite planar graph can be interpreted through the \emph{height} function as a discrete surface~\cite{Thurston}. From this point of view, scaling limits of dimer models on planar bipartite graphs have already been the object of several studies: a law of large number has been established \cite{CKP,KO2} showing that this discrete surface approaches a \emph{limit shape} when the mesh size goes to zero. The fluctuations of the height function around the limit shape have also been studied in the case of a graph embedded in a bounded region, as well as in the case of an \emph{isoradial} infinite graph with critical weights \cite{Bea:biGFF}. The continuous limiting object for these fluctuations is the \emph{Gaussian Free Field}~\cite{Shef:GFF,GlimmJaffe}.  

In this paper, we are interested in the scaling limit of dimer models on $\ZZ^2$-periodic planar graphs but from a different standpoint. Instead of looking at the height function, we consider other observables, called  \emph{pattern density fields}.

Let $G$ be a $\ZZ^2$-periodic planar bipartite graph. A \emph{geometric realization} of $G$ is  an application from $G$ to $\RR^2$ preserving the $\ZZ^2$-periodicity of $G$: vertices of $G$ are mapped to points of $\RR^2$, edges to segments, and $\ZZ^2$ acts on the image of $G$ by translation.

Let $\Psi$ be a geometric realization of $G$ such that the quotient of $\RR^2$ by the action of $\ZZ^2$ has area 1.
For each scaling factor $\epsi>0$, we define the scaled geometric realization $\Psi^\epsi=\epsi\Psi$ and $G^\epsi=\Psi^\epsi(G)$ the image of $G$ by the application $\Psi^\epsi$.


A \emph{pattern} $\pat$ is a finite set of edges $\{\e{e}_1,\ldots,\e{e}_k\}$, together with a marked vertex. The position of the pattern in $G^\epsi$ is given by the coordinates of the image by $\Psi^\epsi$ of this marked vertex.
If $\e{e}_j$ goes from white vertex $\wv_j$ to black vertex $\bv_j$, then the probability of seeing such a pattern in a random dimer configuration is given by formula \eqref{eq:correledges}
\begin{equation*}
 \prob{\pat}= \left(\prod_{j=1}^k \K(\wv_j,\bv_j)\right)\det_{1\leq i,j \leq k}\bigl[\Ki(\bv_i,\wv_j)\bigr].
\end{equation*}

In order to get some information about spatial distribution of patterns, and the way they interact with each other, we define for every pattern $\pat$ a family of (discrete) random fields $N^\epsi_\pat$, called \emph{pattern density fields}. For a given $\epsi$, $N^\epsi_\pat$ is a random distribution, associating to every domain $D$ the number of copies of $\pat$ seen in $D$ in a random dimer configuration of $G^{\epsi}$.
More precisely, if $u^\epsi$ is the image by $\Psi^\epsi$ of the marked vertex attached to $\pat$, the action of $N^\epsi_\pat$ on a smooth test function $\varphi$ is given by
 
\begin{equation}
 N^\epsi_\pat(\varphi)=\sum_{\xx\in\ZZ^2} \varphi(u_{\xx}^\epsi)\mathbbm{I}_{\pat_{\xx}}=\sum_{\xx\in\ZZ^2} \varphi(u_{\xx}^\epsi){\pat_{\xx}}
\end{equation}
where $\pat_{\xx}$ and $u^\epsi_{\xx}$ are the
translates by $\xx$ of $\pat$ and $u^\epsi$. To simplify notations, $\pat_\xx$ will also represent  the indicator function $\mathbbm{I}_{\pat_{\xx}}$ of the pattern $\pat_\xx$, equal to 1 or 0 whether  $\pat_{\xx}$ is in the random dimer configuration $\conf$ or not.

$\epsi^2 \mathbbm{E}\left[N^\epsi_\pat(\varphi)\right]$ is a  Riemann sum of $\mathbbm{P}[\pat]\varphi$, and thus converges to $\mathbbm{P}[\pat]\int\varphi(u)|\ud u|$ when $\epsi$ goes to 0.

The aim of the paper is to prove the following convergence of the fluctuations of this field around its mean value:
\begin{theorem}
When $\epsi$ goes to zero, the normalized fluctuation field
\begin{equation*}
\tilde{N}^\epsi_\pat=\epsi\left(N^\epsi_\pat-\esp{N^\epsi_\pat}\right)
\end{equation*}
converges to a Gaussian field, which is the linear combination of a derivative of the Gaussian Free Field and an independent white noise.
\end{theorem}
See section \ref{sec:statement} for a more precise statement.
Such central limit theorems exist for linear statistics of a broad class of determinantal random fields~\cite{Sosh2} with Hermitian kernels. However, since the kernel $\Ki$ appearing in dimer models is not Hermitian, those results do not apply.
 
Before stating precisely the main results, we recall some basics facts we will need about Gaussian fields.

\subsection{Gaussian Fields}
A Gaussian field is somehow a infinite dimensional generalization of the notion of Gaussian vector. See~\cite{Gelf} for an introduction. As in the classical situation with a Gaussian vector, all the moments can be expressed in terms of the second moment, by  Wick's formula:
\begin{proposition}[Wick's formula]
Let $W$ be a Gaussian field. All the moments of $W$ are determined by the covariance:
Let $\varphi_1,\dots,\varphi_n$ be smooth test functions. Then
\begin{equation}\label{eq:Wick1}
\esp{W(\varphi_1)\cdots W(\varphi_n)}=\begin{cases}
0 & \text{if $n$ is odd,} \\
\displaystyle\sum_{\tau\text{ pairings}} \prod_{k=1}^{n/2} \esp{W(\varphi_{i_k})W(\varphi_{j_k})} & \text{if $n$ is even.}
                                     \end{cases}
\end{equation}
\end{proposition}

In particular, if all the test functions are taken to be equal and $\esp{W(\varphi)^2}=\sigma^2$, then we recover the usual formula for the moments of a Gaussian random variable:
\begin{equation}\label{eq:Wick2}
\esp{W(\varphi)^n}=\begin{cases}
0 & \text{if $n$ is odd,} \\
(n-1)!!\ \sigma^n & \text{if $n$ is even.}
\end{cases}
\end{equation}
where $(n-1)!!=(n-1)\cdots 3\cdot 1$ is the number of pairings of $n$ elements.

%
To prove the convergence in distribution of a sequence of random fields $(W_n)$, in the weak sense, to a Gaussian random field $W$, one has to prove that for all smooth test functions $\varphi$, $W_n(\varphi)$ converges in distribution to $W(\varphi)$. Since the distribution of $W(\varphi)$ is determined by its moments, one has just to check that the moments of $W_n(\varphi)$ converge to that of $W(\varphi)$, \ie that Wick's formula is asymptotically satisfied by $(W_n)$.

\subsection{Statement of the result and outline of the paper\label{sec:statement}}

The main result of the paper is the following Central Limit Theorem when the measure on dimer configurations is fluid (i.e liquid or gaseous).


\begin{theorem}\label{thm:pattern}
Consider the dimer model on a planar bipartite $\ZZ^2$-periodic graph $G$ with a generic liquid or gaseous Gibbs measure $\mu$. Let $\pat$ be a pattern of $G$ and $\flucN{\pat}$ be the random field of density fluctuations of pattern $\pat$. Then when $\epsi$ goes to 0,
$\flucN{\pat}$ converges weakly in distribution to a Gaussian random field $\lflucN{\pat}$.
\begin{itemize}
\item If the measure is liquid, then $\lflucN{\pat}$ is a linear combination of a directional derivative of the massless free field and an independent white noise and its covariance structure has the following form
\begin{equation}
\esp{\lflucN{\pat}(\varphi_1)\lflucN{\pat}(\varphi_2)}= \alpha \iint_{\RR^2\times\RR^2}\!\!\!\!\!\!\!\!\!\!\! \partial\varphi_1(u) \partial\varphi_2(v)G(u,v)\ud u\ud v + \beta \int_{\RR^2}\!\!\!\!\!\varphi_1(u)\varphi_2(u) \ud u 
\end{equation}
where $G(u,v)=-\frac{1}{2\pi}\log|u-v|$ is the Green function on the plane.
\item If $\mu$ is gaseous, then $\lflucN{\pat}$ is a white noise and
\begin{equation}
\esp{\lflucN{\pat}(\varphi_1)\lflucN{\pat}(\varphi_2)}= \beta \int_{\RR^2}\varphi_1(u)\varphi_2(u) \ud u 
\end{equation}
\end{itemize}
\end{theorem}

In other terms, for any choice of
$\varphi_1,\dots,\varphi_n\in\mathcal{C}^\infty_c(\RR^2)$,
the distribution of the random vector
$(\flucN{\pat}\varphi_1,\dots,\flucN{\pat}\varphi_n)$ converges to
that of the Gaussian vector
$(\lflucN{\pat}\varphi_1,\dots,\lflucN{\pat}\varphi_n)$ whose
covariance structure is mentioned in the theorem. As the
distribution of a Gaussian vector is characterized by its moments,
it is sufficient to prove the convergence of the moments of
$(\flucN{\pat}\varphi_1,\dots,\flucN{\pat}\varphi_n)$ to those of
$(\lflucN{\pat}\varphi_1,\dots,\lflucN{\pat}\varphi_n)$, given by Wick's formula~\eqref{eq:Wick1}.

The proof goes in two steps: first we prove the convergence of the second moment of the fluctuation field $\flucN{\pat}$ to the covariance of $\lflucN{\pat}$, and then we prove that Wick's formula is satisfied asymptotically. 

As the arguments are different for a liquid and a gaseous measure, the proof of theorem~\ref{thm:pattern} is decomposed into three cases. In section~\ref{sec:liqedge}, we give the proof for a pattern made of a single edge in the generic liquid case (when the two zeros of $P$ on the unit torus are distinct) and discuss briefly what happens in the non generic case. In section~\ref{sec:liqpat}, the proof is extended to any admissible pattern (\ie a pattern appearing with positive probability) in the generic liquid case. The situation when the measure is  gaseous is discussed in section~\ref{sec:gaz}. The correlations between different pattern density fields are presented in section~\ref{sec:corrfields}, providing a generalization of theorem~\ref{thm:pattern}. In the last section are given two explicit computations : one on the square lattice, and the other on the so-called square-octagon graph. Before entering into the details of the proof, we present in section~\ref{sec:fluiddim} some properties of the different fluid phases of the dimer model.

\section{The fluid phases of the dimer model\label{sec:fluiddim}}

We present here some properties of liquid and gaseous phases of the dimer model on a planar bipartite periodic graph $G$. See~\cite{KOS} for more details. 

\subsection{Liquid phase of dimer models}
\subsubsection*{Map from $G^*$ to $\mathbbm{R}^2$.} 
When given a liquid measure on dimer configurations of $G$, there is a natural
application $\Psi^*$ from the dual graph $G^*$ to $\RR^2$, described by the following lemma. This application seems to give the good
geometry to study liquid dimer models. In particular, when dimer weights are \emph{critical}, $\Psi^*$ coincides with the \emph{isoradial embedding} of $G^*$ described in~\cite{Ke:Crit}.

In the liquid phase, the characteristic polynomial $P$ has two zeros on the unit torus, that are complex conjugate of each other, and generically distinct. Let $(z_0,w_0)$ be one of them.
\begin{lemma}\label{lem:mapping}
 The 1-form 
\begin{equation}
\e{e}=(\wv,\bv)\mapsto i K_{\wv\bv}(z_0,w_0) Q_{\bv\wv}(z_0,w_0)
\end{equation} is a divergence-free flow. Its dual is therefore the gradient of a mapping from $G^*$ to $\RR^2\simeq\CC$.

 This mapping $\Psi^*$ is $\ZZ^2$-periodic and the symmetries of its range are generated by $\hat{\mathrm{x}}=i z_0\partial_1  P(z_0,w_0)$ and $\hat{\mathrm{y}}=i w_0\partial_2  P(z_0,w_0)$.
\end{lemma}

\begin{proof}
The divergence of the 1-form $\omega:\e{e}\mapsto i K_{\wv\bv}(z_0,w_0) Q_{\bv\wv}(z_0,w_0)$ at some black vertex $\bv$ is given by
\begin{align*}
(\div\omega)(\bv)&=\sum_{\wv'\sim\bv}i K_{\wv'\bv}(z_0,w_0) Q_{\bv\wv'}(z_0,w_0)\\&=i\bigl(Q(z_0,w_0)\cdot
K(z_0,w_0)\bigr)_{\bv,\bv}=iP(z_0,w_0)=0.
\end{align*}
Similarly, one can check that the divergence of this flow is 0 at every white vertex $\wv$. Thus, since $G$ is planar, there exists an application $\Psi^*:G^*\rightarrow \CC$ such that $\omega=\ud \Psi$.

The fact that $\Psi^*$ is periodic is a consequence of the fact that $G$ and the Kasteleyn operator are both periodic. There exist two complex numbers $\hat{\mathrm{x}}$ and $\hat{\mathrm{y}}$ such that for every $f\in G^*$ and every $(x,y)\in\ZZ^2$, the difference between the image of $f_{x,y}$ and that of $f$ itself is given by
\begin{equation}
\Psi^*(f_{x,y})-\Psi^*(f)=x\hat{\mathrm{x}}+y\hat{\mathrm{y}}.
\end{equation}
Define $\alpha$ and $\beta$ to be the partial derivatives of $P$ at the root $(z_0,w_0)$ with respect to the first and the second variable respectively
\begin{equation}
\alpha=\partial_1 P(z_0,w_0),\qquad \beta=\partial_2 P(z_0,w_0).
\end{equation}
Let us prove now that the numbers $\hat{\mathrm{x}}$ and $\hat{\mathrm{y}}$ are given by $iz_0\alpha$ and $i w_0\beta$ respectively. Let $\gamma_x$ and $\gamma_y$  be the path delimiting respectively the lower horizontal boundary and the leftmost vertical boundary of the fundamental domain $G_1$ (figure~\ref{fig:crossgamma}).
\begin{figure}
\begin{center}
\includegraphics[width=10cm]{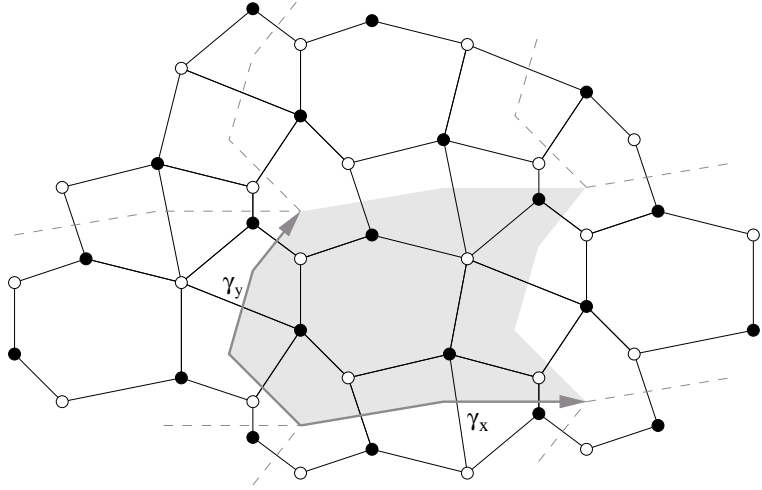}
\caption{\footnotesize A piece of a planar bipartite periodic graph $G$ . The shaded zone represents its fundamental domain $G_1$, delimited by the two paths $\gamma_x$ and $\gamma_y$.}
\label{fig:crossgamma}
\end{center}
\end{figure}
The complex number $\hat{\mathrm{x}}$ equals  the sum of the complex numbers $\pm \omega(\e{e})$ over all edges $\e{e}$ crossing $\gamma_x$ (see figure \ref{fig:crossgamma}). These are exactly the edges of the fundamental domain whose weights have been multiplied by $z^{\pm 1}$ in the Fourier transform of the Kasteleyn operator $K(z,w)$, and the sign of the power of $z$ is the same as that in front of $\omega(\e{e})$. Noticing that
\begin{equation}
z\frac{\partial}{\partial z} z^{m}= m\ z^m=\begin{cases} z & \text{if $m=1$} \\ -z^{-1} & \text{if $m=-1$}\\ 0 & \text{if $m=0$}\end{cases}
\end{equation}
we can write
\begin{align*}
\hat{\mathrm{x}}&=\sum_{\wv,\bv\in G_1} i \bigl(z_0{\partial_1 K_{\wv\bv}}(z_0,w_0)\bigr)Q_{\bv\wv}(z_0,w_0)\\
&=iz_0 \tr\big({\partial_1 K}(z_0,w_0)\cdot Q(z_0,w_0)\bigr).
\end{align*}
On the other hand, since the characteristic polynomial $P(z,w)$ is the determinant of $K(z,w)$, then 
\begin{equation}
i z_0 \alpha= i z_0 \partial_1 P(z_0,w_0)=iz_0 \sum_{j=1}^n \det(K_{\bv_1}(z_0,w_0),\dots,\partial_1 K_{\bv_j}(z_0,w_0),\dots,K_{\bv_n}(z_0,w_0))
\end{equation}
where $K_{\bv_j}(z,w)$ is the $j$th column of the matrix $K(z,w)$. Expanding each determinant with respect to the column containing derivatives, we get
\begin{multline}
i z_0 \alpha= iz_0 \sum_{j=1}^n\sum_{k=1}^n \partial_1 K_{\wv_k\bv_j}(z_0,w_0)\times \\
 \Cof_{\wv_k\bv_j}\bigl(K_{\bv_1}(z_0,w_0),\dots,\partial_1 K_{\bv_j}(z_0,w_0),\dots,K_{\bv_n}(z_0,w_0)\bigr).
\end{multline}
Since the cofactor $\Cof_{\wv_k\bv_j}(K^B_{\bv_1},\dots,\partial_1 K_{\bv_j}(z_0,w_0),\dots,K_{\bv_n}(z_0,w_0))$ does not depend on the $j$th column, we can replace it by $K_{\bv_j}(z_0,w_0)$. This cofactor by definition is $Q_{\bv_j\wv_k}(z_0,w_0)$. Thus,
\begin{align}
i z_0 \alpha&= iz_0\sum_{j=1}^n\sum_{k=1}^n \partial_1 K_{\wv_k\bv_j}(z_0,w_0)Q_{\bv_j\wv_k}(z_0,w_0)\nonumber\\
&=iz_0 \tr(\partial_1 K(z_0,w_0)\cdot Q(z_0,w_0))=\hat{\mathrm{x}}.
\end{align}
The same argument applied to $\gamma_y$ gives the formula
\begin{equation}
i w_0\beta=\hat{\mathrm{y}}.
\end{equation}\qed
\end{proof}
%

In what follows, to construct the application from $G^*$ to $\RR^2$, we will choose between the root $(z_0,w_0)$ on the unit torus for which the frame $(\hat{\mathrm{x}},\hat{\mathrm{y}})$ is direct. This is equivalent to requiring that
\begin{equation}
\label{eq:directframe}
\Im\left(\frac{  w_0  \beta}{z_0  \alpha}\right) >0. 
\end{equation}

To get a geometric realization of $G$ from $\Psi^*$, just pick a point in each dual face of $\Psi^*(G^*)$ in a periodic way. Note that sometimes, as in the case of isoradial embeddings~\cite{Ke:Crit} the obtained realization is not a plane, but can be thought as a globally flat manifold with conic singularities, obtained by gluing together the images of faces of $G^*$.

\subsubsection*{Asymptotics of $\Ki(\bv,\wv)$.}
The coefficients of $\Ki$ decay linearly. More precisely, if $\bv$ and $\wv$ are in the same fundamental domain, and $\bv_{x,y}$ is a translate of $\bv$ by $(x,y)$, then we have the following asymptotics for $\Ki(\bv_{x,y},\wv)$:
\begin{lemma}\label{lem:Kiasymp1}
Let $(z_0,w_0)$ the root of $P$ on the unit torus satisfying \eqref{eq:directframe}. Then the asymptotic expression for the coefficients of $\Ki$ are given by
\begin{multline}\label{eq:Kiasymp1}
\Ki(\bv_{x,y},\wv)= 
-\Re\left(\frac{{z_0}^{-y} {w_0}^{x}Q_{\bv,\wv}(z_0,w_0)}{i\pi\bigl(\alpha z_0 x+\beta w_0 y\bigr)}\right) + O\left(\frac{1}{|x|^2+|y|^2}\right)\\
=-\Re\left(\frac{z_0^{-y} w_0^{x}Q_{\bv,\wv}(z_0,w_0)}{\pi\bigl(x\hat{\mathrm{x}}+y\hat{\mathrm{y}}\bigr)}\right) + O\left(\frac{1}{|x|^2+|y|^2}\right)
\end{multline}
\end{lemma}
The proof of this lemma is given in \cite{KOS}.

Note that the geometry of the map from $G^*$ to $\RR^2$ appears in this analytical result. The denominator up to a factor $\pi$ is the vector separating the image of the fundamental domains of $\bv_{x,y}$ and $\wv$ under this mapping. Moreover, if $\bv$ and $\wv$ are the ends of an edge $\e{e}$ then \eqref{eq:Kiasymp1} can be rewritten as
\begin{equation}
\K(\wv,\bv)\ \Ki(\bv_{x,y},\wv)=
\Re\left(\frac{z_0^{-y} w_0^{x} i\e{e}^*}{\pi\bigl(x\hat{\mathrm{x}}+y\hat{\mathrm{y}}\bigr)}\right) + O\left(\frac{1}{|x|^2+|y|^2}\right)
\end{equation}
and thus correlations between edges decay polynomially with the distance.

\subsection{Gaseous phase}

In this phase, the characteristic polynomial $P$ has no zeros on the unit torus. If $\bv$ and $\wv$ are a black and a white vertex of $G_1$, the fraction ${Q_{\bv\wv}}/{P}$ is analytic on the unit torus and its Fourier coefficients $\Ki(\bv_{\xx},\wv)$ decay exponentially with $\xx$:
\begin{equation*}
\exists\ C_1,C_2>0,\quad\forall\ \xx\in\ZZ^2,\quad\forall\ \bv,\wv\in G_1,\quad \bigl|\Ki(\bv_{\xx},\wv)\bigr| \leq C_1 e^{-C_2 |\xx|}.
\end{equation*}

Hence the correlations also decay exponentially.

\subsection{Density fields and partition function}

In this paragraph, the calculations are purely formal and try to give some heuristics on the information one could get from these pattern density fields.

 Consider a dimer model for which we assign to a configuration
  $\mathcal{C}$ the weight  $w(\mathcal{C})$. Let $Z_0$ be the partition function of this model
\begin{equation}
Z_0=\sum_{\mathcal{C}} w(\mathcal{C})
\end{equation}
We now perturb the partition function by modifying locally the configuration weights. Let $\varphi$ be a smooth test function. Fix a pattern $\pat$ and a $\epsi>0$. We multiply every weight $w(\mathcal{C})$ by a factor $e^{t\epsi\varphi(u^\epsi_\xx)}$ whenever there is a copy of a given pattern $\pat$ located at $\xx$. Up to a multiplicative constant $\exp(t\epsi \prob{\pat}\sum_{\xx}\varphi(u^\epsi_\xx))$, the new partition function for the model with these new weights is
\begin{equation*}
Z_t=\sum_{\mathcal{C}} w(\mathcal{C}) e^{t\epsi \sum_{\xx} \varphi(u^\epsi_\xx)(\mathbbm{I}_{\pat_\xx}-\prob{\pat})}=\sum_{\mathcal{C}}w(\mathcal{C}) e^{t\flucN{\pat}(\varphi)}
\end{equation*}
This can be generalized to a perturbation involving several patterns.
Formally, the successive derivatives of $Z_t/Z_0$ at $t=0$ are the moments of the random variable $\flucN{\pat}(\varphi)$ with respect to the unperturbed probability measure.
\begin{equation*}
\left.\frac{\ud^k }{\ud t^k} \frac{Z_t}{Z_0}\right|_{t=0}=\espind{0}{(\flucN{\pat}\ \varphi)^k}
\end{equation*}
Thus, an information on the moments of $\flucN{\pat}$ gives an infinitesimal information on the perturbed partition function that could be hopefully integrated to construct probability measures corresponding to the modified weights.


\section{The liquid case: edge densities} \label{sec:liqedge}

We first concentrate on the proof of a simple particular case of theorem \ref{thm:pattern}.
We suppose that the probability measure on dimer configurations of $G$ is a generic liquid measure, and consider the fluctuations of the density random field $\flucN{\e{e}}$ for a pattern consisting in a single edge $\e{e}=(\wv,\bv)$.
The precise statement we prove in this section is the following:
\begin{theorem}\label{thm:liqedge}
The random field $\flucN{\e{e}}$ converges weakly in distribution, as $\epsi$ goes to 0, to a Gaussian random field $\lflucN{\e{e}}$ with covariance
\begin{equation}
\esp{\lflucN{\e{e}}(\varphi_1)\lflucN{\e{e}}(\varphi_2)}=\frac{1}{\pi}\iint \partial_{\e{e}^*}\varphi_1(u) G(u,v)\partial_{\e{e}^*}\varphi_2(v) \ud u\ud v+ A\int \varphi_1(u)\varphi_2(u)\ud u
\end{equation}
for a certain $A\geq 0$, and where $\e{e}^*$ is the vector representing the dual edge of $\e{e}$ in the geometric realization $\Psi^*$ of $G^*$ presented in lemma~\ref{lem:mapping}. 
\end{theorem}
Note that the particular geometry introduced by the application $\Psi$ is particularly well-adapted to the problem.

As we discussed in the previous section, we first prove the convergence of the second moment, and then that of higher moments. 
In this section, as we are interested in copies of an
 edge $\e{e}=(\wv,\bv)$, the only vertices we will deal with are most of the time translates of $\wv$ and $\bv$. To simplify
 notations, we will write $\Ki(\xx-\xx')$ instead of
 $\Ki(\bv_{\xx},\wv_{\xx'})$ and $\K_{\e{e}}$ will stand for $\K(\wv,\bv)$.

\subsection{Convergence of the second moment\label{subsec:cov}}

The second moment $(\varphi_1,\varphi_2)\mapsto\esp{\flucN{\e{e}}(\varphi_1)\flucN{\e{e}}(\varphi_2)}$
of $\flucN{\e{e}}$ is a continuous bilinear positive form on
 $\mathcal{C}^\infty_c(\RR^2)$.
We prove that this bilinear form converges to a non
degenerate bilinear form, that will define the covariance structure
for the limit Gaussian field $\lflucN{\e{e}}$.

\begin{proposition}
\label{prop:cov}
There exists a non-negative constant $A$ such that
\begin{multline}
 \forall \varphi_1,\varphi_2 \in \mathcal{C}^\infty_c(\RR^2)
\\ \lim_{\epsi\rightarrow 0}
\esp{\flucN{\e{e}}(\varphi_1)\flucN{\e{e}}(\varphi_2)}=
\frac{1}{\pi}\iint_{\RR^2\times\RR^2 }\partial_{\e{e}^*}\varphi_1(u_1)
G(u_1,u_2) \partial_{\e{e}^*} \varphi_2(u_2)|\ud u_1||\ud u_2|\\ +A
\int_{\RR^2} \varphi_1(u)\varphi_2(u_2) |\ud u|
\end{multline}
where $G(u,v)=-\frac{1}{2\pi}\log|u-v|$ is the Green function on the plane.
\end{proposition}

Before going into the proof of this proposition, we give some interpretation of the expression for the covariance.
The right hand side  can be physically interpreted as the energy of interaction between
two magnetic dipoles with moment density $\varphi_1\e{e}^*$ and
$\varphi_2\e{e}^*$, plus a term of interaction at very short range.

Suppose that there is an excess of edges $\e{e}$ in the random dimer configuration, in some region $D_1$. 
These dimers behave collectively as a magnetic dipole: their presence influences the rest of the dimer configuration as if a magnetic field created by a dipole with a density $\e{e}^*$ over $D_1$ was applied to the system: edges whose dual has an orientation closed to that of the magnetic field at that point is more likely to appear.  

\begin{proof}
Using the invariance by translation of the Kasteleyn operator and
hence of the correlations, we rewrite the second moment as a
convolution of two distributions, applied to a test function.

\begin{align}
\esp{\flucN{\e{e}}(\varphi_1)\flucN{\e{e}}(\varphi_2)}&=\epsi^2\sum_{\xx,
\xx'\in
\ZZ^2}\varphi_1(u_{\xx})\varphi_2(u_{\xx'})\esp{(\e{e}_{\xx}-\overline
{\e{e}})(\e{e}_{\xx'}-\overline{\e{e}})} \\
&=\epsi^2 \sum_{\xx,\xx'\in
\ZZ^2}\varphi_1(u_{\xx})\varphi_2(u_{\xx'})\esp{(\e{e}
-\overline{\e{e}})(\e{e}_{\xx-\xx'}-\overline{\e{e}})}  \\
 &= \langle \varphi_1^\epsi \ast F^\epsi,\varphi_2\rangle
\end{align}
where $\overline{\e{e}}=\prob{\e{e}}$ and the two distributions $\varphi_1^\epsi$ and $F^\epsi$ are
defined by

\begin{equation*} \varphi_1^\epsi=\epsi^2\sum_{xy}\varphi_1(u_{\xx})
\delta_{u_{\xx}^\epsi} \qquad
F^\epsi=\sum_{\xx}\esp{(\e{e}-\overline{\e{e}})(\e{e}_{\xx}
-\overline{\e{e}})}\delta_{u_{\xx}^\epsi}.
\end{equation*}

The distribution $\varphi_1^\epsi$ converges weakly to $\varphi_1$ when $\epsi$ goes
to zero. We will now prove the convergence of $F^\epsi$ to some distribution
$F$, what will ensure that $\varphi_1^\epsi \ast F^\epsi$ converges
weakly to $\varphi_1^\epsi \ast F^\epsi$, since the support of
$\varphi_1^\epsi$ is contained in the fixed compact
$\mathrm{supp}(\varphi_1)$, and hence that $\langle \varphi_1^\epsi \ast F^\epsi,\varphi_2\rangle$ converges.

Let $\psi$ be a smooth test function with compact support. Let us prove the convergence of

\begin{equation}\label{eq:Feps}
\langle F^\epsi,\psi\rangle=\sum_{\xx}\esp{(\e{e}-\overline{\e{e}})(\e{e}_{\xx}
-\overline{\e{e}})}\psi(u_{\xx}^\epsi)=\sum_{\xx}\mathrm{Cov}(\e{e},\e{e}_{\xx}
)\psi(u_{\xx}^\epsi).
\end{equation}

At first sight, $\langle F^\epsi, \psi\rangle$ looks vaguely like a
Riemann sum of a particular function. The problem is that, due to the
asymptotics of $\Ki$, the function would behave as $1/u^2$, which is
not integrable in the vicinity of 0. Therefore, we decompose the sum
on $\xx$ in the definition of $F^\epsi$ depending on whether
the norm of $\xx$ is larger than $M=\lfloor 1/ \epsi\rfloor$
or not, that is if $u^\epsi_{\xx}$ is in $\mathcal{B}=\left\{-i(\alpha z_0 s +\beta w_0 t) \ ;\
(s,t)\in[-1,1]^2\right\}$ or not.
\begin{multline} \label{eq:Feps3sums}
\langle F^\epsi,\psi\rangle = \sum_{|\xx|>M}
\psi(u^\epsi_{\xx})\mathrm{Cov}(\e{e},\e{e}_{\xx})
+ \sum_{|\xx|\leq M}
\psi(u^\epsi_{\xx})\mathrm{Cov}(\e{e},\e{e}_{\xx}) \\= \sum_{|\xx|>M}
\psi(u^\epsi_{\xx})\mathrm{Cov}(\e{e},\e{e}_{\xx})
+\sum_{|\xx|\leq M}
\left(\psi(u^\epsi_{\xx})-\psi(0)\right)\mathrm{Cov}(\e{e},\e{e}_{\xx}) +\psi(0)\sum_{|\xx|\leq M}
\mathrm{Cov}(\e{e},\e{e}_{\xx}).
\end{multline}

The fact we subtracted and added $\psi(0)$ in the second sum
removed the non integrable singularity at 0.
The following two lemmas state the convergence of the three sums.

\begin{lemma}
\begin{align}
\lim_{\epsi\rightarrow 0}\sum_{|\xx|>\lfloor1/\epsi\rfloor}
\psi(u^\epsi_{\xx})\mathrm{Cov}(\e{e},\e{e}_{\xx})
&=\frac{1}{2\pi^2}\int_{\RR^2\setminus
\mathcal{B}}\psi(z)\Re{\left(\frac{\e{e}^*}{u}\right)^2}|\ud u| \label{eq:limsum1}\\
\lim_{\epsi\rightarrow 0}\sum_{|\xx|\leq \lfloor1/\epsi\rfloor}
\left(\psi(u^\epsi_{\xx})-\psi(0)\right)\mathrm{Cov}(\e{e},\e{e}_{\xx})&=
\frac{1}{2\pi^2}\int_{\mathcal{B}}(\psi(z)-\psi(0))\Re\left(\frac{\e{e}^*} { u}\right)^2|\ud u| \label{eq:limsum2}
\end{align}
Moreover the sum of the two previous limits can be rewritten as
\begin{equation*}
\frac{1}{\pi} \int_{\RR^2}\partial_{\e{e}^*}\psi(u)\partial_{\e{e}^*}G(u,0)|\ud u| + A_1\psi(0).
\end{equation*}
for some constant $A_1$.
\end{lemma}

%
\begin{proof}
If $\xx\neq(0,0)$, the covariance between edges $\e{e}=\e{e}_{(0,0)}$ and $\e{e}_{\xx}=\e{e}_{(x,y)}$
is given by
\begin{align*}
\mathrm{Cov}(\e{e},\e{e}_{\xx})&=\prob{\e{e}\textrm{ and }\e{e}_\xx}-\prob{\e{e}}^2 =\K_{\e{e}}^2 \det\begin{bmatrix} \Ki(0) & \Ki(\xx) \\ \Ki(-\xx) & \Ki(0)\end{bmatrix}-\left(\K_{\e{e}} \Ki(0)\right)^2
\\&=-\K_{\e{e}}^2\ \Ki(\xx)\Ki(-\xx)
\end{align*}

Using asymptotics of $\Ki$ for large $\xx$, we get the following asymptotic
expression for the covariance between two distinct edges
\begin{align*}
\mathrm{Cov}(\e{e},\e{e}_{\xx})&=-\frac{\epsi^2}{\pi^2}\Re\left(\frac{z_0^{-y}w_0^{x}i\e{e}^*}{u_\xx^{\epsi}}
\right)\Re\left(\frac{z_0^{y}w_0^{-x}i\e{e}^*}{-u_{\xx}^{\epsi}}
\right)+O\left(\frac{\epsi}{u^{\epsi}_{\xx}}\right)^3 \\
&=-\frac{\epsi^2}{2\pi^2}\Re\left(\frac{\e{e}^*}{
u^{\epsi}_{\xx}}\right)^2+\frac{\epsi^2}{2\pi^2}\Re(z_0^{2y}w_0^{-2x}
)\left|\frac{\e{e}^*}{u^{\epsi}_{\xx}}\right|^2+O\left(\frac{
\epsi}{u^{\epsi}_{\xx}}\right)^3
\end{align*}

Since the second term is oscillating, it will not contribute to the limit. The
sum in the left hand side of \eqref{eq:limsum1}, modulo the oscillating terms, can be interpreted as the integral of a
piecewise constant function approximating
\begin{equation*}
-\frac{1}{2\pi^2}\psi(u)\Re\left(\frac{{\e{e}^*}^2}{u^2}\right)\mathbbm{I}_{
\RR^2\setminus \mathcal{B}}(u)
\end{equation*}

As the approximating functions are bounded uniformly in $\epsi$ by an
integrable
function, and converge almost everywhere, then by
Lebesgue theorem, the first sum converges to
\begin{equation}
-\int_{\RR^2\setminus
\mathcal{B}}\frac{1}{2\pi^2}\Re{\left(\frac{{\e{e}^*}^2}{u^2}\right)}\psi(z)|\ud
u|.
\end{equation}

In the same way, the sum in the left hand side of \eqref{eq:limsum2} is the integral of a piecewise constant function
approximating
\begin{equation*}
-\frac{1}{2\pi^2}(\psi(u)-\psi(0))\Re\left(\frac{{\e{e}^*}^2}{
u^2}\right)\mathbbm{I}_{\mathcal{B}}(u)
\end{equation*}
and for the same reasons, it converges to
\begin{equation}
-\int_{\mathcal{B}}\frac{1}{2\pi^2}\Re\left(\frac{{\e{e}^*}^2}{u^2}
\right)\bigl(\psi(z)-\psi(0)\bigr)|\ud u|.
\end{equation}

We rewrite the sum of the limit using Green's formula inside and outside of $\mathcal{B}$, noticing that
$\frac{1}{2\pi}\Re\left(\frac{\e{e}^*}{u}\right)^2$ is the second derivative of the Green function $G(u,0)=-\frac{1}{2\pi}\log|u|$ along the vector $\e{e}^*$.
\begin{multline}
\int_{\RR^2\setminus
\mathcal{B}}\frac{1}{\pi}\partial^2_{\e{e}^*}G(u,0)\psi(z)|\ud
u| =\\
\frac{1}{\pi}\oint_{\partial\mathcal{B}}\psi(u)\partial_{\e{e}^*}G(u,0)\langle \vec{n}_{in},\e{e}^*\rangle \ud\sigma
-\frac{1}{\pi}\int_{\RR^2\setminus\mathcal{B}}\partial_{\e{e}^*}\psi(u)
\partial_{\e{e}^*}G(u,0)|\ud u|,
\end{multline}
\begin{multline}
\int_{\mathcal{B}}\frac{1}{\pi}\partial^2_{\e{e}^*}G(u,0)(\psi(u)-\psi(0))|\ud u|=
\\ \frac{1}{\pi}\oint_{\partial\mathcal{B}}(\psi(u)-\psi(0))\partial_{\e{e}^*}G(u,0)\langle \vec{n}_{out},\e{e}^*\rangle \ud\sigma
-\frac{1}{\pi}\int_{\mathcal{B}}\partial_{\e{e}^*}\psi(u) \partial_{\e{e}^*}G(u,0)|\ud u|.
\end{multline}
where $\vec{n}_{in}$ and $\vec{n}_{out}$ are the unit normal vector fields on $\partial \mathcal{B}$ pointing respectively inwards and outwards. The two
integrals on $\mathcal{B}$ and $\RR^2\setminus\mathcal{B}$
combine to give an integral over $\RR^2$. The two integrals
on $\partial \mathcal{B}$ cancel out partially. It remains
only
\begin{equation*}
-\psi(0)\oint_{\partial\mathcal{B}}\partial_{\e{e}^*}G(u,0)\langle \vec{n}_{ext},\e{e}^*\rangle \ud\sigma = A_1 \psi(0).
\end{equation*}
This establishes the convergence of the two first sums , which completes the proof of the lemma.\qed
\end{proof}

We now prove the convergence of the third sum.
\begin{lemma}\label{lem:convfour}
$\sum_{|\xx|\leq M} \mathrm{Cov}(\e{e},\e{e}_{\xx})$ converges when
$M$ goes to infinity to  a limit $A_2$.
\end{lemma}
\begin{proof}
The sum of the covariances is given in terms of $\K$ and $\Ki$ by
\begin{multline}
 \sum_{|\xx|\leq M}  \mathrm{Cov}(\e{e},\e{e}_{\xx})
=\prob{\e{e}}\left(1-\prob{\e{e}}\right)+\sum_{0<|\xx|\leq M} \mathrm{Cov}(\e{e},\e{e}_{\xx})
\\=\K_{\e{e}}\Ki(0)-\K_{\e{e}}^2\sum_{|\xx|\leq M}.
\Ki(\xx)\Ki(-\xx)
\end{multline}

 $\Ki(\xx)$ is  by definition the
$\xx$th Fourier coefficient of the function $f=Q/P$ defined on the unit torus $\mathbbm{T}^2$.
As $P$ has simple zeros, $f$ is in $\mathrm{L}^1(\mathbbm{T}^2)$. The
convolution
\begin{equation}
(f\ast f)(z,w)=\iint_{\mathbbm{T}^2} f(u,v)f({z}{u},{w}{v})\frac{\ud
u}{2 i\pi u}\frac{\ud v}{2 i\pi v}
\end{equation}
is also in $\mathrm{L}^1(\mathbbm{T}^2)$ and its $\xx$th Fourier coefficient is
exactly $\Ki(\xx)\Ki(-\xx)$. Establishing the convergence of the
sum is now a problem of pointwise convergence of a Fourier series. If $f$ had been
continuous at $(z,w)=(1,1)$, then the Fourier series would have converged to
$f(1,1)$. The problem is that the function
$f(u,v)^2$ is not integrable and thus, the function $f\ast f$ is not
defined when $z$ and $w$ are both equal to 1. However, $f\ast f$ is smooth in a punctured
neighborhood of $(1,1)$, has directional limits when $(z,w)$ converges to $(1,1)$, varying
continuously with the direction. We can then prove an analogue in two dimensions
of Dirichlet's theorem
for $f\ast f$ to show that the Fourier series at $(z,w)=(1,1)$ when $M$ goes to infinity, converges to a mixture of the directional limits of $f\ast f$. More precisely, if $t=\arg(w)/\arg(z)$ and $\ell(t)$ the limit of
$f\ast f$ when $(z,w)$ goes to $(1,1)$ with $t$ fixed, then

\begin{equation}
\lim_{M\rightarrow +\infty}\sum_{|\xx|\leq M}
\Ki(\xx)\Ki(-\xx) =
\frac{1}{\pi^2}\int_{-\infty}^{+\infty} \ell(t)\log\left|\frac{1+t}{1-t}\right|
\frac{\ud t}{t}
\end{equation}
And thus $\sum_{|\xx|\leq M} \mathrm{Cov}{(\e{e}-\overline{\e{e}})(\e{e}_{\xx}-\overline{\e{e}})}$ converges to a limit that we denote by $A_2$.\qed
\end{proof}

We now come back to the proof of the convergence of the distribution $F^\epsi$. The three sums in \eqref{eq:Feps3sums}
defining $\langle F^\epsi,\psi\rangle$ converge and the sum of the limits is
\begin{equation}
\langle F,\psi\rangle = \frac{1}{\pi} \int_{\RR^2}\partial_{\e{e}^*}\psi(u)\partial_{\e{e}^*}G(u,0)|\ud u| + A\psi(0)
\end{equation}
where $A=A_1+A_2$. Thus, when $\epsi$ goes to 0,
$F^\epsi$ converges to the distribution $F$ defined by the
formula above, and hence $\varphi_1^\epsi\ast F^\epsi$
to $\varphi_1 \ast F$. Denoting by
$\partial^{\scriptscriptstyle(u)}$ and
$\partial^{\scriptscriptstyle(v)}$ respectively the operator of
partial differentiation with respect to the variable $u$ and
$v$, and noticing that, since $G(u,v)=G(u-v)$, we have:
\begin{equation*}
\partial^{\scriptscriptstyle(u)}G(u,v)= -\partial^{\scriptscriptstyle(v)}G(u,v),
\end{equation*}
\noindent
we get the following expression for the limit covariance structure
\begin{align}
\lim_{\epsi \rightarrow 0} &\esp{(\flucN{\e{e}}\varphi_1)(\flucN{\e{e}}\varphi_2)}=\langle \varphi_1 \ast F,\varphi_2\rangle \nonumber\\
&=\frac{1}{\pi}\iint_{\RR^2\times\RR^2} \partial_{\e{e}^*}\varphi_1(u) \partial_{\e{e}^*}^{\scriptscriptstyle(u)}G(u,v) \varphi_2(v)|\ud u||\ud v|+A \int_{\RR^2} \varphi_1(u)\varphi_2(u)|\ud u|\nonumber\\
&=-\frac{1}{\pi}\iint_{\RR^2\times\RR^2} \partial_{\e{e}^*}\varphi_1(u) \partial_{\e{e}^*}^{\scriptscriptstyle(v)}G(u,v) \varphi_2(v)|\ud u||\ud v|+A \int_{\RR^2} \varphi_1(u)\varphi_2(u)|\ud u|\nonumber\\
&=\frac{1}{\pi}\iint_{\RR^2\times\RR^2} \partial_{\e{e}^*}\varphi_1(u) G(u,v) \partial_{\e{e}^*}\varphi_2(v)|\ud u||\ud v|+A \int_{\RR^2} \varphi_1(u)\varphi_2(u)|\ud u|
\end{align}

Thus the covariance of $\flucN{\e{e}}$ converges to the expression given in proposition \ref{prop:cov}.\qed
\end{proof}


\subsection{Convergence of higher moments}\label{subsec:highmmt}

We now prove the convergence of the moments of order $\geq 3$ of $\flucN{\e{e}}$ to those of the Gaussian field $\lflucN{\e{e}}$.
\begin{proposition} \label{prop:nth1} For every $n\geq 3$,
The $n$th moment of $\flucN{\e{e}}$ converges to that of $\lflucN{\e{e}}$ when $\epsi$ goes to zero. In other words,
for every $\varphi_1,\cdots,\varphi_n \in \mathcal{C}^\infty_c(\RR^2)$,
\begin{equation*}
\lim_{\epsi\rightarrow 0}
\esp{\flucN{\e{e}}{\varphi_1}\cdots\flucN{\e{e}}{\varphi_n}}=\esp{\lflucN{\e{e}}{\varphi_1}\cdots\lflucN{\e{e}}{\varphi_n}}
\end{equation*}
\end{proposition}
Since $\lflucN{\e{e}}$ is Gaussian, it is sufficient to show that in
the limit, the moments of $\flucN{\e{e}}$ satisfy Wick's formula.
Moreover, as
$(\varphi_1,\cdots,\varphi_n)\mapsto\esp{\flucN{\e{e}}{\varphi_1}\cdots\flucN{\e
{e}}{\varphi_n}}$ is a symmetric $n$-linear form, we just have to
prove proposition \ref{prop:nth1} when all the $\varphi_i$ are equal to some test
function $\psi$, the general case being obtained by polarization.
The previous proposition reduces then to showing that

\begin{proposition}\label{prop:nth}
\begin{equation*}
\lim_{\epsi\rightarrow 0} \esp{(\flucN{\e{e}}\psi)^n}=\esp{(\lflucN{\e{e}}\psi)^n}=
\left\{\begin{array}{cl} 0 &\textrm{if $n$ is odd} \\
(n-1)!!  \ \esp{(\lflucN{\e{e}}\psi)^2}^{ n/2} &\textrm{if $n$ is even}
\end{array}\right.
\end{equation*}

\end{proposition}

In this section, we are thus interested in the
limit of

\begin{equation}
\label{eq:nthmmt}
\esp{(\flucN{\e{e}}\psi)^n}
=\epsi^n \sum_{\xx_1,\cdots,\xx_n}
\psi(u^\epsi_1)\cdots\psi(u^\epsi_n)\esp{(\e{e}_{\xx_1}-\overline{\e
{e}})\cdots(\e{e}_{\xx_n}-\overline{\e{e}})}
\end{equation}

A first step in the proof is to study the convergence of a related quantity $\Xi_n^\epsi(\psi)$, defined by a sum of the same general term  as for $\esp{(\flucN{\e{e}}\psi)^n}$, but with a set of indices $(\xx_j)$ restricted to distinct points:

\begin{equation}\label{eq:Xi_n}
\Xi_n^\epsi(\psi)=\epsi^n\sum_{\substack{\xx_1\cdots \xx_n \\ \textrm{distinct}}} \psi(u^\epsi_{\xx_1})\cdots\psi(u^\epsi_{\xx_n})\esp{(\e{e}_{\xx_1}-\overline{\e{e}})\cdots(\e{e}_{\xx_n}-\overline{\e{e}})}.
\end{equation}



\subsubsection{Convergence of $\Xi_n^\epsi(\psi)$}

To prove the convergence of $\Xi_n^\epsi(\psi)$, we have to understand the asymptotic behavior of the correlations between distinct edges, when they are far from each other. A simple expression is given by Kenyon in \cite{Ke:ConfInv}
 to compute these correlations using a unique
determinant.

\begin{lemma} [\cite{Ke:ConfInv}]
\label{lem:corr}
Let $\e{e}_1=(\wv_1,\bv_1),\cdots,\e{e}_n=(\wv_n,\bv_n)$ be distinct edges. Their correlation is given by
\begin{align*}
\esp{(\e{e}_1-\bar{\e{e}}_1)\cdots(\e{e}
_n-{\bar{\e{e}}_n})}&=\left(\prod_{j=1}^n
\K(\wv_j,\bv_j)\right)\det_{1\leq i,j \leq n}\begin{bmatrix} 0 & &
\Ki(\bv_i,\wv_j) \\ & \ddots & \\ \Ki(\bv_j,\wv_i) & & 0 \end{bmatrix}
\end{align*}
\end{lemma}


This formula allows us to give an explicit expression for $\Xi_n^\epsi(\psi)$
in terms of the operators $\K$ and $\Ki$. Since the matrix in lemma \ref{lem:corr} has zeros on the diagonal, only permutations with no fixed point will contribute to the determinant expressed as a sum over the symmetric group. Let $\hat{\mathfrak{S}}_n$ be the set of such permutations. Every permutation $\sigma\in\hat{\mathfrak{S}}_n$ is decomposed as a product of disjoint cycles $\gamma_1\cdots\gamma_p$. The supports of these cycles form a partition $(\Gamma_l)_{l=1}^p$ of $\{1,\dots,n\}$, whose parts $\Gamma_l$ have cardinal at least 2. The terms coming from permutations leading to the same partition are put together and we get:
\begin{multline}\label{eq:dvpmmt}
\epsi^n\sum_{\substack{\xx_1\cdots \xx_n \\ \textrm{distinct}}} \esp{\prod_{j=1}^n \psi(u^\epsi_{\xx_j}) (\e{e}_{\xx_j}-\bar{\e{e}})}
=\epsi^n\K_{e}^n\sum_{\substack{\xx_1\cdots \xx_n \\ \textrm{distinct}}} \prod_{j=1}^n \psi(u^\epsi_{\xx_j})\det\begin{bmatrix} \scriptscriptstyle 0 & 
\scriptscriptstyle\Ki(\bv_i,\wv_j) \\ \scriptscriptstyle\Ki(\bv_j,\wv_i)  & \scriptscriptstyle 0 \end{bmatrix}\\
= \epsi^n\K_{e}^n\sum_{\substack{\xx_1\cdots \xx_n \\ \textrm{distinct}}}\sum_{\sigma\in
\hat{\mathfrak{S}}_n} \sgn(\sigma)\prod_{j=1}^n \psi(u^\epsi_{\xx_j}) \Ki({\xx_{\sigma(j)}}-{\xx_j}) \\
\!\!\!\!\!\!\!=\sum_{\{\Gamma\}_{l=1}^p}\prod_{l=1}^p\sum_{\substack{\gamma\textrm{
cycle}\\
\mathrm{supp}(\gamma)=\Gamma_l}}\left(\sgn(\gamma)\epsi^{|\gamma|}\sum_{\substack{\xx_{j_1},\cdots,\xx_{j_{|\gamma|}}
\\ \text{distinct}}}\prod_{k=1}^{|\gamma|}
\psi(u^\epsi_{\xx_{j_k}})\K_{\e{e}}\Ki({\xx_{\gamma(j_k)}}-{\xx_{j_k}})
+o(1)\right)\!\!\!\!
\end{multline}

 The error term $o(1)$ comes from the fact that in the second line, we allow two $\xx_j$ whose indices are in different components of $(\Gamma_l)$ to be equal.

 We now examine the convergence of a term in brackets, associated to a cycle $\gamma$
\begin{equation*}
\sgn(\gamma)\epsi^{|\gamma|}\sum_{\substack{\xx_{j_1},\cdots,\xx_{j_{|\gamma|}}
\\ \text{distinct}}}\prod_{k=1}^{|\gamma|}
\psi(u^\epsi_{\xx_{j_k}})\K_{\e{e}}\Ki({\xx_{\gamma(j_k)}}-{\xx_{j_k}})
\end{equation*}
According to subsection \ref{subsec:cov}, we know that when $\gamma$ is a transposition the corresponding term converges. When the length of $\gamma$ is at least 3, we have the following lemma.

\begin{lemma}
For any  cycle $\gamma$ of length $m\geq 3$, and any $\psi\in\mathcal{C}^{\infty}_c(\RR^2)$, we have
\begin{multline}
\lim_{\epsi\rightarrow 0}\epsi^m \sum_{\substack{\xx_{1},\cdots \xx_m \\ \textrm{distinct}}}
\prod_{k=1}^m\psi(u^\epsi_{\xx_k})\K_{\e{e}}
\Ki(\xx_{\gamma(k)}-\xx_{k})=\\
\frac{1}{2^{m-1}\pi^m}\int_{(\RR^2)^m}
\Re\left(\frac{(i\e{e}^*)^m\psi(u_1)\cdots\psi(u_m)}{(u_{\gamma(1)}-u_{1})\cdots(u_{\gamma(m)}-u_{m})}\right)|\ud u_1|\cdots|\ud u_m|
\end{multline}
\end{lemma}

\begin{proof}
From the behavior of $\Ki$ at long range, we get an asymptotic
expression for the product
\begin{align} \label{eq:Xiedgeasymp}
\prod_{j=1}^m\K_{\e{e}}\Ki(\xx_{\gamma(j)}-\xx_j)&=\prod_{j=1}^m\Re\left(\frac{\epsi z_0^{-(y_{\gamma(j)}-y_j)}w_0^{+(x_{\gamma(j)}-x_j)}i\e{e}^*} {\pi(u_{\xx_{\gamma(j)}}-u^\epsi_{\xx_{j}})}\right)+\smallterms
\\ &=\frac{\epsi^m}{2^{m-1}\pi^m
}\Re\left(\frac{(i\e{e}^*)^m}{(u^\epsi_{\xx_{\gamma(1)}}-u^\epsi_{\xx_{1}})\cdots(u^\epsi_{\xx_{\gamma(m)}}-u^\epsi_{\xx_{m}})}\right) +
\oscterms
\end{align}

The oscillating part of this asymptotic expansion once summed over $\xx_1,\dots,\xx_m$ will not contribute to the limit.
The sum of the leading term multiplied by $\epsi^m\psi(u_{\xx_1}^\epsi)\cdots\psi(u_{\xx_m}^\epsi)$ can be interpreted as the integral of a piecewise constant function, converging almost everywhere to

\begin{equation}
\frac{1}{2^{m-1}\pi^m}\psi(u_1)\cdots\psi(u_m)
\Re\left(\prod_{j=1}^m\frac{i\e{e}^*}{u_{\gamma(j)}-u_j}\right)
\end{equation}

As all the functions are dominated by a constant times the integrable function
\begin{equation*}
\frac{|\psi(u_1)\cdots\psi(u_m)|}{|(u_{\gamma(1)}-u_1)\cdots (u_{\gamma(m)}-u_m)|},
\end{equation*}
 the convergence follows from Lebesgue theorem.\qed
\end{proof}

Once proven the convergence of all these terms, we can combine their limit to get the limit of $\Xi_n^\epsi(\psi)$. When summing over all cycles with a given support $\Gamma=\{j_1,\dots,j_m\}$ of cardinality $m\geq 3$, we get the following limit :

\begin{multline*}
\lim_{\epsi\rightarrow 0}\sum_{\substack{\gamma\textrm{ cycle}\\ \mathrm{supp}(\gamma)=\Gamma}}\left(\epsi^{m}\sum_{\substack{\xx_{j_1},\cdots,\xx_{j_{m}} \\distinct}}\prod_{k=1}^{m} \psi(u^\epsi_{\xx_{j_k}})\K_{\e{e}}\Ki(\xx_{\gamma(j_k)}-\xx_{j_k}) +o(1)\right)\\
=\frac{1}{2^{m-1}\pi^m}\int_{(\RR^2)^{m}}\prod_{j\in\Gamma}\psi(u_j)\Re\left(\sum_{\substack{\gamma\textrm{ cycle}\\ \mathrm{supp{\gamma}=\Gamma}}}
\prod_{j\in\Gamma}\frac{1}{u_{\gamma(j)}-u_{j}}\right) |\ud u_1|\cdots|\ud u_m|
\end{multline*}

which equals zero according to the following lemma :

\begin{lemma}
Let $m\geq 3$, and $u_1,\ldots,u_n$ be distinct complex numbers. Then
\begin{equation}\label{eq:lemcycles}
\sum_{\substack{\gamma\in\mathfrak{S}_m \\ m-cycles}}\prod_{i=1}^{m}\frac{1}{u_{\gamma(i)}-u_{i}}=0
\end{equation}
\end{lemma}

\begin{proof}
Denote by $f$ the function of $u_1,\dots,u_m$ defined by the left hand side of \eqref{eq:lemcycles}.
When $m$ is odd, the result is obvious, since $\gamma$ and $\gamma^{-1}$ give opposite contributions.
For a general $m$, since $m$-cycles form a conjugation class in the group of permutations $\mathfrak{S}_m$, the function $f$ is a rational fraction invariant under permutation of the variables $u_1,\dots,u_m$. The denominator of this fraction is the Vandermonde $V=\prod_{i<j}(u_i-u_j)$, and the numerator is of lower degree than $V$. Since $V$ is antisymmetric under permutation, the denominator has to be as well. But the only antisymmetric polynomial of lower degree than $V$ is 0.\qed
\end{proof}

In the limit of equation \eqref{eq:dvpmmt} will contribute only partitions whose all components have cardinality 2. If $n$ is odd, $\{1,\dots,n\}$ cannot be partionned into parts of two elements, and the limit $\Xi_n(\psi)$ of $\Xi_n^\epsi(\psi)$ is zero. When $n$ is even,
there are $(n-1)!!$ such partitions, corresponding each to a pairing. The limit $\Xi_n(\psi)$ of $\Xi_n^\epsi(\psi)$ is then

\begin{align*}
\Xi_n(\psi)&=\sum_{\substack{\textrm{pairings} \\ \{i_1,j_1\},\dots,\{i_{n/2},j_{n/2}\}}} \prod_{k=1}^{n/2} \left(\lim_{\epsi\rightarrow 0} \epsi^2 \sum_{\xx_{i_k}\neq \xx_{j_k}}\psi(u^\epsi_{\xx_{i_k}})\psi(u^\epsi_{\xx_{j_k}})\mathrm{Cov}(\e{e}_{\xx_{i_k}},\e{e}_{\xx_{j_k}})\right)\\
&=(n-1)!! \ \left( \Xi_2(\psi)\right)^{n/2}
\end{align*}

 For the moment, we discussed the limit of a restricted sum on distinct edges, which is not exactly the expression for the $n$th moment. We will now deal with the case when some edges can coincide.

\subsubsection{Proof of proposition \ref{prop:nth}}
\label{ssec:nth2}
In the expression of the $n$th moment
\begin{equation*}
\epsi^n \sum_{\xx_1,\cdots,\xx_n}
\psi(u^\epsi_1)\cdots\psi(u^\epsi_n)\esp{(\e{e}_{\xx_1}-\overline{\e
{e}})\cdots(\e{e}_{\xx_n}-\overline{\e{e}})}
\end{equation*}
the correlations are not given by lemma \ref{lem:corr} as soon as at least two edges coincide.
To understand the behavior of this expression, we must be able to express in terms of $\K$ and $\Ki$ correlations of the form
\begin{equation}
\label{eq:gencorr}
\esp{(\e{e}_{\xx_1}-\overline{\e{e}})^{k_1}\cdots(\e{e}_{\xx_p}-\overline{\e{e}}
)^{k_p}}
\end{equation}
when $\e{e}_1,\dots,\e{e}_p$ are distinct, and $k_1,\dots,k_p\geq 1$.
Using the fact that the indicator function of an edge $\e{e}$ satisfies
$\e{e}^j=\e{e}$ for $j\geq 1$, Newton formula yields

\begin{equation*}
(\e{e}-\overline{\e{e}})^k=\sum_{j=0}^k \binom{k}{j} \e{e}^j
(-\overline{\e{e}})^{k-j}=\e{e}\sum_{j=1}^k \binom{k}{j}
(-\overline{\e{e}})^{k-j} +(-\overline{\e{e}})^k=\alpha_k^{\e{e}}
(\e{e}-\overline{\e{e}})+ \beta_k^{\e{e}}
\end{equation*}
where $\alpha_k^{\e{e}}$ and $\beta_k^{\e{e}}$ are deterministic, depending only
on
$\overline{\e{e}}$ and $k$. Since all the edges we consider are translates one from another, and thus have the same
probability, we will simply denote these coefficients by $\alpha_k$ and
$\beta_k$. Note some particular values of $\alpha_k$ and $\beta_k$ that will be useful later
\begin{equation*}
 \alpha_1=1, \qquad \beta_1=0, \qquad
\beta_2=\overline{\e{e}}(1-\overline{\e{e}}).
\end{equation*}

Correlation \eqref{eq:gencorr} can be rewritten with these notations as
\begin{align}
\esp{(\e{e}_{\xx_1}-\overline{\e{e}})^{k_1}\cdots(\e{e}_{\xx_p}-\overline{\e{e}}
)^{k_p}}&=\esp{(\alpha_{k_1}(\e{e}_1-\overline{\e{e}})+\beta_{k_1}
)\cdots(\alpha_{k_p}(\e{e}_p-\overline{\e{e}})+\beta_{k_p})} \nonumber \\
&=\sum_{J\subset\{1,\cdots,p\}} \tilde{\alpha}_J \tilde{\beta}_{\bar{J}}
\esp{\prod_{j\in J} (\e{e}_j-\overline{\e{e}})}
\end{align}
where $\displaystyle \tilde{\alpha}_J=\prod_{j\in J}\alpha_{k_j}$
and$\displaystyle \tilde{\beta}_{\bar{J}}=\prod_{j\notin J}\beta_{k_j}$.

In equation \eqref{eq:nthmmt}, $\esp{(\flucN{\e{e}}\psi)^n}$ is expressed as
a sum over all edges. We want now to rewrite this sum as a sum over distinct
edges, using partitions of $\{1,\dots,n\}$. Such a partition $\{\Delta_l\}_{l=1}^p$ is associated naturally to every $n$-tuple of lattice points $(\xx_1,\dots,\xx_n)$:  each component $\Delta_l$ of this partition is an
equivalence class for the relation
\begin{equation*} i \sim j \Leftrightarrow \xx_i = \xx_j.
\end{equation*}
 Denoting by $n_l$ the cardinal of $\Delta_l$, we rewrite  equation \eqref{eq:nthmmt} as
\begin{multline}
\label{eq:nthmmt2}
\esp{(\flucN{\e{e}}\psi)^n} =\epsi^n\sum_{\{\Delta_l\}_{l=1}^p}
\sum_{\substack{\xx_1,\cdots,\xx_p\\ \textrm{distinct}}}  \esp{\prod_{j=1}^p\psi^{n_j}(u^\epsi_{\xx_j})(\e{e}_{\xx_j}-\bar{\e{e}})^{n_j}} \\
=\epsi^n\sum_{\{\Delta_l\}_{l=1}^p}
\sum_{J\subset\{1,\dots,p\}}\tilde{\alpha}_J\tilde{\beta}_{\bar J}
\sum_{\substack{(\xx_j)_{j\in J}\\ \textrm{distinct}}} \esp{\prod_{j\in J}
\psi^{n_j}(u^\epsi_{\xx_j})(\e{e}_{\xx_j}-\bar{\e{e}})}
{\sum_{(\xx_l)_{l\notin J}}}' \prod_{l\notin J}
\psi^{n_l}(u^\epsi_{\xx_l})
\end{multline}
where $\Sigma'$ means the sum over $(\xx_l)$ distinct, but also distinct from values of any $\xx_j,j\in J$.
Denote by $q$ the number of $\Delta_l$ reduced to a single element. As $\beta_1=0$, $\tilde{\beta}_{\bar J}$ is zero unless $J$ contains the indices of these $\Delta_l$. Thus, the cardinal  of a subset $J$ giving a non-zero contribution must be at least q.  For such a $J$, the last sum over $(\xx_l)_{l\notin J}$ in \eqref{eq:nthmmt2} is a
Riemann sum, and therefore is $O(\epsi^{-2(p-|J|)})$.


Furthermore, the sum over $(\xx_J)_{j\in J}$ can be expressed by polarization in terms of $\Xi_{|J|}^\epsi$, and is therefore a $O(\epsi^{-|J|})$.
Since $|J|\geq q$ and
$$n=\sum_{l=1}^p n_l\geq q+2(p-q)=2p-q,$$
 the contribution of $J$ to \eqref{eq:nthmmt2} is at most
$O(\epsi^{n-2p+|J|})$ which will be negligible in the limit except when $|J|=q$ and $n_l=2$ for all $l\notin J$. For such $J$ and $(\Delta_l)$, we have
\begin{equation}
\tilde{\alpha}_J=1,\qquad \tilde{\beta}_{\bar J}=\left(\bar{\e{e}}(1-\bar{\e{e}})\right)^{(n-q)/2}.
\end{equation}
Thus the only partitions that will contribute to the limit are ``partial pairings'', matching by pairs $2m=(n-q)$ elements  of $\{1,\dots,n\}$. For a fixed $m$, there are
$$\binom{n}{2m}(2m-1)!!=\frac{n!}{2^m m! (n-2m)!}$$
 such partitions, all giving the same contribution. Summing over $m$ we get
\begin{multline*}
 \esp{(\flucN{\e{e}}\psi)^n}=\sum_{m=0}^{\lfloor n/2\rfloor}\binom{n}{2m}(2m-1)!!\  \Xi^\epsi_{n-2m}(\psi)  \left(\epsi^2\bar{\e{e}}(1-\bar{\e{e}})\sum_{\xx}\psi^2(u^\epsi_{\xx})\right)^m+O(\epsi).
\end{multline*}

The Riemann sum $\epsi^2\bar{\e{e}}(1-\bar{\e{e}})\sum \psi^2(u_\xx^\epsi)$ converges to
\begin{equation*}
\bar{\e{e}}(1-\bar{\e{e}})\int \psi^2(u)\ud u=\esp{\lflucN{\e{e}}(\psi)^2}-\Xi_2(\psi).
\end{equation*}
If $n$ is odd, so is $n-2m$. In this case, $\lim\Xi^\epsi_{n-2m}(\psi) =0$, and therefore
\begin{equation}
\lim_{\epsi\rightarrow 0}\esp{(\flucN{\e{e}}\psi)^n}=0.
\end{equation}

If $n$ is even,
\begin{align}
\lim_{\epsi\rightarrow 0}\Xi^\epsi_{n-2m}(\psi) &=\Xi_{n-2m}(\psi)=(n-2m-1)!!\ \Xi_2(\psi)^{n/2-m}\\
&=\frac{(n-2m)!}{2^{n/2-m}(n/2-m)!}\Xi_2(\psi)^{n/2-m}.
\end{align}

Therefore, the limit of $\esp{(\flucN{\e{e}}\psi)^n}$ is given by
\begin{align*}
\lim_{\epsi\rightarrow 0}\esp{(\flucN{\e{e}}\psi)^n}&=\sum_{m=0}^{ n/2} \frac{n!}{2^m m! (n-2m)!}\frac{(n-2m)!\ \Xi_2(\psi)^{n/2-m}}{2^{n/2-m}(n/2-m)!} \left(\bar{\e{e}}(1-\bar{\e{e}})\int\psi^2(u)|\ud u|\right)^m\\
&=(n-1)!!\sum_{m=0}^{ n/2}\binom{n/2}{m}\ \Xi_2(\psi)^{n/2-m}\left(\bar{\e{e}}(1-\bar{\e{e}})\int\psi^2(u)|\ud u|\right)^m \\
&=(n-1)!!\left(\Xi_2(\psi)+\bar{\e{e}}(1-\bar{\e{e}})\int\psi^2(u)|\ud u|\right)^{n/2}
\\&=(n-1)!!\left(\esp{\lflucN{\e{e}}(\psi)^2}\right)^{n/2}
\end{align*}
what is exactly what we wanted to prove. This therefore ends the proof of theorem \ref{thm:pattern} for a pattern made of one edge and a generic liquid Gibbs measure.

\subsection{A word about the non generic case}

When the two roots of the characteristic polynomial $P(z,w)$ on the unit torus coincide, the measure is still liquid, and the correlations between edges at distance $r$ still decay like $r^{-2}$. However, since $z_0$ and $w_0$ are real, the leading term is the asymptotics of $\Ki$ is not oscillating anymore, what will induce a ``resonance phenomenon'' in the system.

The two first sums in equation \eqref{eq:Feps3sums} defining the distribution $F^\epsi$ appearing in the study of the convergence of the second moment still have a finite limit when $\epsi$ goes to zero. On the contrary, due to the resonance, the third sum
\begin{equation*}
\sum_{|\xx|\leq \lfloor 1/\epsi \rfloor} \mathrm{Cov}(\e{e},\e{e}_{\xx})
\end{equation*}
in this case diverges. More precisely, this sum is $O(\log(1/\epsi))$. Therefore the second moment diverges. However, one can prove that $(\log(1/\epsi))^{-1}\flucN{\e{e}}$ converges weakly in distribution to a white noise. We will not show it here.

\section{The liquid case: general pattern densities}\label{sec:liqpat}

This section is devoted to the proof of an analogue of theorem \ref{thm:liqedge} for multi-edged patterns. Let $\flucN{\pat}$ the density fluctuation field of a pattern $\pat$ for a generic liquid Gibbs measure.

\begin{theorem}
When $\epsi$ goes to zero, $\flucN{\pat}$ converges weakly in distribution to the Gaussian field $\lflucN{\pat}$ whose covariant structure is given by
\begin{multline*}
\esp{\lflucN{\pat}(\varphi_1)\lflucN{\pat}(\varphi_2)}=\frac{1}{\pi}\iint_{\RR^2\times\RR^2}\partial_{\pat^*}\varphi_1(u_1)G(u_1,u_2) \partial_{\pat^*}\varphi_2(u_2)|\ud u_1||\ud u_2|\\
+A\int_{\RR^2}\varphi_1(u_1)\varphi_2(u)|\ud u|
\end{multline*}
where the vector $\pat^*$ and the nonnegative constant $A$ depend only on the Gibbs measure and the pattern $\pat$.
\end{theorem}

The scheme of the proof is very similar to that of theorem \ref{thm:liqedge} for edge density fluctuations. The problem is that there is no simple analogue of lemma \ref{lem:corr} for correlations between non overlapping patterns, that is why the proof needs a little more combinatorial work in this case.

After having introduced the different notations required to deal easily with these patterns, we prove the theorem, following the structure of the proof given in the last section, and explaining in details only parts that are specific to patterns made of more than one edge.

\subsection{Notations}

Let $\pat$ be a pattern containing $k$ distinct edges $\e{e}^1=(\wv^1,\bv^1),\dots,\e{e}^k=(\wv^k,\bv^k)$. The probability of such a pattern to appear in a random dimer configuration is

\begin{equation}
\overline{\pat}=\prob{\pat}= \K_{\e{e}^1}\cdots \K_{\e{e}^k} \det \left[ \begin{array}{ccc} \Ki(\bv^1,\wv^1) & \cdots & \Ki(\bv^1,\wv^k) \\ \vdots & \ddots&\vdots \\\Ki(\bv^k,\wv^1) & \cdots & \Ki(\bv^k,\wv^k) \end{array} \right].
\end{equation}

More generally, the probability to see $n$ non-overlapping copies $\pat_1,\dots,\pat_n$ of $\pat$ obtained respectively by translation of a lattice vector $\xx_1,\dots,\xx_n$ is given up to a constant by a determinant of matrix $nk\times nk$ defined by blocks

\begin{equation}
\prob{\pat_1\cdots\pat_n}= \left(\K_{\e{e}^1}\cdots \K_{\e{e}^k}\right)^n \det\left[\begin{array}{c|c|c} A_{11} &\cdots &A_{1n} \\\hline \vdots& \ddots&\vdots\\ \hline A_{n1} & \cdots&A_{nn} \end{array}\right]
\end{equation}
where the entries of the block $A_{IJ}$ are coefficients of $\Ki$ between black vertices of $\pat_I$ and white vertices of $\pat_J$

\begin{equation}
A_{IJ}=\left[ \begin{array}{ccc} \Ki(\bv^1_{\xx_I},\wv^1_{\xx_J}) & \cdots & \Ki(\bv^1_{\xx_I},\wv^k_{\xx_J}) \\ \vdots &\ddots &\vdots \\\Ki(\bv^k_{\xx_I},\wv^1_{\xx_J}) & \cdots & \Ki(\bv^k_{\xx_I},\wv^k_{\xx_J}) \end{array} \right].
\end{equation}

The matrix $A_{II}$ does not depend on $I$. We denote by $E$ this matrix, whose determinant is used to compute $\bar{\pat}=\prob{\pat}$. We suppose that the pattern appears with positive probability $\bar{\pat}>0$, and in particular, $E$ is invertible. Defining  $B_{IJ}$ as the product $E^{-1} A_{IJ}$ and $B$ as the whole block matrix $(B_{IJ})$, we can rewrite the joint probability of $\pat_1,\dots,\pat_n$ as

\begin{equation} \label{eq:blockcorrel}
\prob{\pat_1\cdots\pat_n}= (\overline{\pat})^n \det\left[\begin{array}{c|c|c|c} \mathbbm{I}_k &B_{1,2}&\cdots &B_{1,n} \\\hline B_{2,1} & \mathbbm{I}_k & \phantom{B_{11}}& \vdots\\ \hline \vdots&\phantom{B_{11}} &\ddots & \vdots\\ \hline B_{n,1} &\cdots  & B_{n,n-1}&\mathbbm{I}_{k} \end{array}\right].
\end{equation}

Instead of using a single integer $i$ to denote the line (resp. the column) of an entry in such a matrix defined by blocks, it will be more convenient to use a couple of integers $(I,\alpha)$, where $I$ is the index of the block line (resp. of the block column) and $\alpha$ is the relative position in the $I$th block line (resp. block column). The relation between the two sets of indices is simply $$i=I(k-1)+\alpha.$$

If the coordinates $\xx_j$ are all distinct but some patterns partially overlap, define
$\tilde{\pat}_j=\pat_j\setminus{\bigcup_{i=1}^{j-1}\pat_i}$ for all
$j\in\{1,\dots,n\} $. We have then
\begin{equation*}
\prob{\pat_1\cdots\pat_n}=\prob{\tilde{\pat}_1\cdots\tilde{\pat}_n}.
\end{equation*}
Up to a relabeling of the patterns, we can assume that none of the
$\tilde{\pat}_j$ is empty. Thus the joint probability of these patterns is also given by the determinant of a matrix defined by blocks of size $|\tilde{\pat}_1|+\cdots+|\tilde{\pat}_n|$.

\subsection{Asymptotics of correlations}

The following lemma gives asymptotic correlations between distant patterns.

\begin{lemma}\label{lem:blockasymp}
Let $(\xx_j)=\xx_1,\dots,\xx_n$ be $n$ distinct lattice points. The correlations between the patterns $\pat_{\xx_1},\dots,\pat_{\xx_n}$ can be rewritten as
\begin{equation*}
\esp{(\pat_{\xx_1}-\bar{\pat})\cdots(\pat_{\xx_n}-\bar{\pat})}=\sum_{S\in\hat{\mathfrak{S}}_n}\left(\prod_{\gamma\text{ cycle of }S} H_{\gamma}\bigl((\xx_j)\bigr)\right)+\smallterms
\end{equation*}
where the functions $H_{\gamma}$ have the following asymptotic behaviour
\begin{equation*}
H_{\gamma}\bigl((\xx_j)\bigr)=\sgn(\gamma)\frac{\epsi^{|\gamma|}}{2^{|\gamma|-1}\pi^{|\gamma|}}\Re\left(\frac{\tr((E^{-1}Q)^{|\gamma|})}{
\prod_{j\in\mathrm{supp}\gamma}u^\epsi_{\xx_{\gamma(j)}}-u^\epsi_{\xx_{j}}}\right)+\oscterms
\end{equation*}
with $Q=\bigl(Q_{\bv^\alpha \wv^\beta}(z_0,w_0)\bigr)_{\alpha,\beta=1}^k$, 
and satisfy
\begin{equation*}
\left|H_{\gamma}\bigl((\xx_j)\bigr)\right|\leq \epsi^{|\gamma|}\frac{C}{\prod_{j\in\mathrm{supp}\gamma}|u_{\gamma(j)}-u_{j}|}
\end{equation*}
for every $u_1,\dots,u_n$ in a $\epsi$-neighborhood of $u^\epsi_{\xx_1},\dots,u^\epsi_{\xx_n}$.
The error term $o(1)$ is uniformly bounded in $\xx_1,\dots,\xx_n$ and goes to zero when the distance between the patterns goes to infinity.

\end{lemma}


\begin{proof}
We first derive the asymptotic expression for the correlations when the patterns are far from each other. When $|\xx_i-\xx_j|$ is large enough for every $i\neq j$, the patterns are disjoint and expression \eqref{eq:blockcorrel} for correlations can be used. Expanding the products in the expectation, we get
\begin{multline}
\esp{(\pat_{\xx_1}-\bar{\pat})\cdots(\pat_{\xx_n}-\bar{\pat})}=\sum_{j_1,\dots,j_p}(-\bar{\pat})^{n-p}\esp{\pat_{\xx_{j_1}}\dots\pat_{\xx_{j_p}}}\\
=\prob{\pat}^n\sum_{C\subset\{1,\dots,n\} } (-1)^{n-|C|}\det \left[\begin{array}{c|c|c|c} \mathbbm{I}_n & \delta_{12} B_{12}& \cdots &\delta_{1n} B_{1n} \\
\hline \delta_{21} B_{21} & \mathbbm{I}_n & & \\
 \hline \vdots& &\ddots &\\
  \hline \delta_{n1}B_{n1} & & & \mathbbm{I}_n\end{array} \right]
\end{multline}
where the non diagonal bloc $(I,J)$ is either $B_{IJ}$ or $0$ depending on whether $(I,J)$ belongs to $C\times C$ or not.

Expressing each determinant as a sum over the symmetric group $\mathfrak{S}_{nk}$ and gathering the terms coming from the same permutation, one can notice that the contributions of permutations fixing a whole block are vanishing, due to the alternating sign in the sum over $C$. The permutations contributing to the correlations are those whose support intersects each block. Therefore, we have
\begin{equation*}
\esp{(\pat_{\xx_1}-\bar{\pat})\cdots(\pat_{\xx_n}-\bar{\pat})}=\prob{\pat}^n\sum_{\substack{\sigma\in\mathfrak{S}_{nk}\\
\textrm{fixing no block}}} \sgn(\sigma)  \prod_{I,\alpha}
B_{(I,\alpha),\sigma((I,\alpha))}
\end{equation*}

The main contribution to this sum is given by the ``special'' permutations, the support of which intersects each block exactly once. 
Let $\sigma$ be such a permutation. The $n$ non fixed elements are $(1,\alpha_1),\dots(n,\alpha_n)$. 
For every $I\in\{1,\dots,n\}$, the non fixed element $(I,\alpha_I)$ is sent to $\sigma((I,\alpha_I))=(S(I),\alpha_{S(I)}$, where $S\in\hat{\mathfrak{S}}_n$ is a fixed-point free permutation, having the same signature as $\sigma$. There are $k^n$ special permutations $\sigma$ leading to the same $S$, corresponding to the different possible choices of the non fixed points in each block $\alpha_1,\dots,\alpha_n$.

The contribution of the other permutations will be negligible because of the extra $\epsi$ coming from additional entries of $\K^{-1}$ between non-fixed points.

\begin{align}
\esp{\prod_{j=1}^n(\pat_{\xx_j}-\bar{\pat})}&=\prob{\pat}^n\sum_{\substack{\sigma\in\mathfrak{S}_{nk}\\
\textrm{fixing no block}}} \sgn(\sigma)  \prod_{I,\alpha}
B_{(I,\alpha),\sigma((I,\alpha))}\nonumber\\
&=\prob{\pat}^n \sum_{S\in\hat{\mathfrak{S}}_n} 
\sgn(S) \left(\prod_{I=1}^n\sum_{\alpha_I=1}^k
(B_{I,S(I)})_{\alpha_I,\alpha_{S(I)}}\right)+\smallterms
\label{eq:blockasymp}
\end{align}

The main term in equation \eqref{eq:blockasymp} has an expression in terms of traces of products of block matrices $B_{IJ}$
\begin{equation}
\prod_{I=1}^n\sum_{\alpha_I=1}^k\left(B_{I,S(I)}\right)_{\alpha_I\alpha_{S(I)}}=\prod_{\substack{\gamma=(I_1,\dots,I_p)\\
\textrm{cycle of }S}}\tr\left(B_{I_1,I_2}\cdots B_{I_p,I_1}\right)
\end{equation}

Let us have a look to a particular trace $\tr\left(B_{I_1,I_2}\cdots
B_{I_p,I_1}\right)$. Recall that $B_{IJ}$ is the product of $E^{-1}$
whose entries will be denoted by $e_{\alpha\beta}$ and $A_{IJ}$
whose entries $(A_{IJ})_{\alpha\beta}$ are  the
coefficient of $\Ki$ taken between black vertex `$\alpha$' of
pattern $I$ and white vertex `$\beta$' of pattern $J$. The
asymptotics of $(A_{IJ})_{\alpha\beta}$ when patterns are far away from each other are given by
lemma \ref{lem:Kiasymp1}:
\begin{equation}
(A_{IJ})_{\alpha\beta}=\Ki(\bv^\alpha_{\xx_I},\wv^\beta_{\xx_J})=-\frac{\epsi}{\pi}\Re\left(\frac{z_0^{-y_J+y_I}w_0^{+x_J-x_I}Q_{\alpha\beta}(z_0,w_0)}{u^\epsi_{\xx_J}-u^\epsi_{\xx_I}}\right)+\smallterms.
\end{equation}

To simplify notations, we introduce
\begin{equation}
\zeta_{IJ}=\frac{z_0^{-y_J+y_I}w_0^{+x_J-x_I}}{u^\epsi_{\xx_J}-u^\epsi_{\xx_I}}
\end{equation}
and write $Q_{\alpha\beta}$ instead of $Q_{\bv^{\alpha}\wv^{\beta}}(z_0,w_0)$.
 The trace $\tr\left(B_{I_1,I_2}\cdots B_{I_p,I_1}\right)$ can therefore be rewritten as
\begin{multline*}
\tr\left(B_{I_1,I_2}\cdots B_{I_p,I_1}\right) = \sum_{\alpha_1,\dots,\alpha_k}(E^{-1}A_{I_1I_2})_{\alpha_1\alpha_2}\cdots(E^{-1}A_{I_p I_1})_{\alpha_p\alpha_1} \\
=\sum_{\substack{\alpha_1,\dots,\alpha_k \\ \beta_1,\cdots,\beta_k}}\left(\frac{-\epsi}{\pi}\right)^p e_{\alpha_1\beta_1}\cdots e_{\alpha_p \beta_p}\Re\left(Q_{\beta_1\alpha_2}\zeta_{I_1 I_2}\right)\cdots \Re\left(Q_{\beta_p\alpha_1}\zeta_{I_p I_1}\right) +\smallterms.
\end{multline*}

As in the case of edge densities, there will be only two non
oscillating terms in the expansion of this product of real parts:
those for which the phases contained in the $\zeta_{IJ}$ compensate
exactly. 

With the convention that $p+1=1$, one has
\begin{multline*}
\tr\left(B_{I_1 I_2}\right. \cdots \left.B_{I_p I_1}\right)= \\= \sum_{\substack{\alpha_1,\dots,\alpha_k \\ \beta_1,\cdots,\beta_k}}\left(\frac{-\epsi}{2\pi}\right)^p\prod_{j=1}^p e_{\alpha_j\beta_j}\left(Q_{\beta_j\alpha_{j+1}}\zeta_{I_j I_{j+1}}+ \overline{Q_{\beta_j\alpha_{j+1}}\zeta_{I_j I_{j+1}}}\right)+ \smallterms \\
=\sum_{\substack{\alpha_1,\dots,\alpha_k \\ \beta_1,\cdots,\beta_k}}\frac{(-\epsi)^p}{2^p\pi^p}\left(\prod_{j=1}^p e_{\alpha_j\beta_j}Q_{\beta_j\alpha_{j+1}}\zeta_{I_j I_{j+1}}+\prod_{j=1}^p e_{\alpha_j\beta_j} \overline{Q_{\beta_j\alpha_{j+1}}\zeta_{I_j I_{j+1}}}\right)+\oscterms.
\end{multline*}

The product of the $\zeta_{I_kI_{k+1}}$ is equal to
\begin{equation}
\prod_{j=1}^p \zeta_{I_j,I_{j+1}}=\prod_{j=1}^p\frac{1}{ u^\epsi_{I_{j+1}}-u^\epsi_{I_{j}}}
\end{equation}
and we can rewrite the trace of $B_{I_1,I_2}\cdots B_{I_p,I_1}$ as

\begin{equation}
\tr\left(B_{I_1,I_2}\cdots B_{I_p,I_1}\right)
=\frac{(-\epsi)^p}{2^{p-1}\pi^p}\Re\frac{\tr((E^{-1}Q)^p)}{(u^\epsi_{I_{2}}-u^\epsi_{I_{1}})\dots(u^\epsi_{I_{1}}-u^\epsi_{I_{p}})}+\oscterms
\end{equation}
giving the asymptotics for $H_{\gamma}$.

When the patterns are not disjoint anymore,  then a similar analysis can be done, in defining new patterns as the connected components of $\pat_1\cup\dots\cup\pat_n$. The bound can be extended to this case.\qed
\end{proof}

%

\subsection{Convergence of the second moment}

\begin{proposition}\label{prop:liqpat2moment}
\begin{multline}
\forall \ \varphi_1,\varphi_2 \in \mathcal{C}^{\infty}_c(\RR^2)\quad \\
\lim_{\epsi \rightarrow 0}\esp{\flucN{\pat}(\varphi_1) \flucN{\pat}(\varphi_2)}
=\frac{1}{\pi} \iint_{\RR^2\times\RR^2} \partial_{\pat^*}\varphi_1(u_1) G(u_1,u_2) \partial_{\pat^*}\varphi_2(u_2) |\ud u_1||\ud u_2| \\+ A \int_{\RR^2} \varphi_1(u)\varphi_2(u) |\ud u|
\end{multline}
\end{proposition}

The proof of the convergence of the second moment goes exactly as
that of section \ref{sec:liqedge}. The second moment of
$\flucN{\pat}$ can be expressed as a convolution of two
distributions applied to a test function

\begin{equation}
\esp{(\flucN{\pat}\varphi_1)(\flucN{\pat}\varphi_2)}=\langle
\varphi^\epsi\ast F^\epsi ,\varphi_2\rangle
\end{equation}

with the same definitions as before for $\varphi_1^{\epsi}$ and
$F^\epsi$.

\begin{equation}
\varphi_1^{\epsi}=\epsi^2\sum_{\xx}\varphi_1(u^{\epsi}_{\xx})
\delta(\cdot-u^{\epsi}_{\xx}) \qquad
F^\epsi=\sum_{\xx}\mathrm{Cov}(\pat,\pat_{\xx})
\delta(\cdot-u^{\epsi}_{\xx})
\end{equation}

$\varphi_1^{\epsi}$ converges weakly to $\varphi_1$. The
convergence of $F^\epsi$ to a distribution $F$ is proven exactly in the same way as in section \ref{sec:liqedge}. The only
difficulty that could appear is the analogue of lemma \ref{lem:convfour}
proving the convergence of
$\sum_{\xx}\mathrm{Cov}(\pat,\pat_{\xx})$.
$\mathrm{Cov}(\pat,\pat_{\xx})$ is a linear combination of products
of diverse values of $\Ki$. If we interpret these products as
Fourier coefficients of a convolution of functions, then the
convergence becomes more obvious: we saw in section
\ref{sec:liqedge} that the function whose Fourier coefficients are
the product of two $\Ki$ is not defined at $(1,1)$, but has
directional limits when $(z,w)$ converges to $(1,1)$, and the sum of
the Fourier coefficient was an average of these directional limits.
When more than two $\Ki$ are involved, the function is even continuous
thanks to the multiple convolutions, and the Fourier series
converges at $(z,w)=1$.

$F^\epsi$ converges then to the distribution $F$ acting on a test function $\psi$ as
\begin{equation}
\langle F, \psi\rangle=\frac{1}{\pi}\iint_{\RR^2} \psi(u)\Re\frac{\tr\bigl((E^{-1}Q)^2\bigr)}{2\pi u^2} \ud u + A \psi(0)
\end{equation}

The complex number representing the vector $\pat^*$ along which are
taken the derivatives is a square root of $\tr(E^{-1}QE^{-1}Q)$, and after application of Green's formula, we get the expression of proposition \ref{prop:liqpat2moment} for $\langle \varphi_1\ast F,\varphi_2\rangle$.

\subsection{Convergence of higher moments}

\begin{proposition}
Let $n\geq 3$, and $\varphi_1,\dots,\varphi_n\in\mathcal{C}^{\infty}_c(\RR^2)$.

\begin{equation}
\lim_{\epsi \rightarrow 0} \esp{\flucN{\pat}(\varphi_1)\cdots \flucN{\pat}(\varphi_n)}=
\left\{\begin{array}{c l} 0 & \textrm{if $n$ is odd} \\
\displaystyle \sum_{\text{pairings}} \prod_{l=1}^{n/2}
\esp{\lflucN{\pat}(\varphi_{i_l})\lflucN{\pat}(\varphi_{j_l})} & \textrm{if $n$ is even} \end{array} \right.
\end{equation}

\end{proposition}

\begin{proof}

As usual, it is sufficient to study the case where all the $\varphi_i$ are equal to some fixed smooth test function $\psi$. The $n$th moment is the given by

\begin{equation}
\esp{(\flucN{\pat}\psi)^n}=\epsi^n\sum_{\xx_1,\dots,\xx_n}\esp{\prod_{l=1}^n\psi(u^{\epsi}_l)(\pat_{\xx_l}-\bar{\pat})}.
\end{equation}
We know from lemma \ref{lem:blockasymp} the asymptotics of the general term of the sum

\begin{align*}
\epsi^n&\esp{\prod_{l=1}^n\psi(u^{\epsi}_l)(\pat_{\xx_l}-\bar{\pat})}\\
&=\sum_{S\in\hat{\mathfrak{S}}_n}\prod_{\gamma \text{ cycle of }S}
\sgn{\gamma}\frac{2\epsi^{2|\gamma|}} {(2\pi)^{|\gamma|}}
\Re\left(\tr\bigl((E^{-1}Q)^{|\gamma|}\bigr)\prod_{l\in\mathrm{supp}\gamma}\frac{\psi(u^{\epsi}_l)}{
u^\epsi_{\xx_{\gamma(l)}}-u^\epsi_{\xx_{l}}}\right) \\
&\qquad\qquad\qquad\qquad\qquad\qquad\qquad\qquad\qquad\qquad\qquad\qquad {+\oscterms+\smallterms}\\
&=\sum_{(\Gamma_l)_1^p}\prod_{l=1}^p  (-2)\left(\frac{-\epsi^2}{2\pi}\right)^{|\Gamma_l|}
\Re\left(\tr\bigl((E^{-1}Q)^{|\Gamma_l|}\bigr)\sum_{ \substack{\gamma\text{ cycle} \\  \mathrm{supp}\gamma=\Gamma_l}}\prod_{j\in\Gamma_l}
 \frac{\psi(u^{\epsi}_j)}{ u^{\epsi}_{\xx_{\gamma(j)}}-u^{\epsi}_{\xx_{j}}}\right) \\
 &\qquad\qquad\qquad\qquad\qquad\qquad\qquad\qquad\qquad\qquad\qquad\qquad{+\smallterms+\oscterms}
\end{align*}
where the $(\Gamma_l)$ are partitions of $\{1,\dots,n\}$ whose components have size at least 2. The expression we obtained is very close to that of equation \eqref{eq:Xiedgeasymp}. From this point, the same arguments as for edges yield the proof of the proposition.\qed \end{proof}



\section{The gaseous case}\label{sec:gaz}

In a gaseous phase, $\Ki$ decays exponentially. There exist two constants $C_1$ and $C_2$ such that

$$\forall \ \xx,\xx'\in\ZZ^2,\quad \left|\Ki(\bv^{i}_{\xx},\wv^{j}_{00})\right|\leq C_1\cdot e^{-C_2 |\xx|}$$

The precise statement of theorem \ref{thm:pattern} in this particular context for a pattern consisting in a single edge $\e{e}$ is the following:
\begin{theorem}\label{thm:gasedge}
The random field $\flucN{\e{e}}$ converges weakly in distribution to a white noise of amplitude $\sqrt{\frac{\partial^2 \mathcal{F}}{\partial^2 \log \K_{\e{e}}}}$, where
\begin{equation}
\mathcal{F}=\iint_{\tore} \log P(z,w)\frac{\ud z}{2i\pi z}\frac{\ud w}{2i\pi w} 
\end{equation}
is the free energy per fundamental domain of the dimer model.
\end{theorem}

The proof is exposed in the first two subsections, and the case of a more complex pattern is briefly discussed in subsection \ref{subsec:gaspat}
\subsection{Convergence of the second moment}

\begin{proposition}
\begin{equation*}
\forall \varphi_1,\varphi_2 \in \mathcal{C}_c^\infty(\RR^2) \quad
\lim_{\epsi\rightarrow
0}\esp{\flucN{\e{e}}(\varphi_1)\flucN{\e{e}}(\varphi_2)}=
\frac{\partial^2 \mathcal{F}}{\partial \log \K_{\e{e}}^2} \int_{\RR^2}
\varphi_1(z)\varphi_2(z) |\ud z|.
\end{equation*}
\end{proposition}

\begin{proof}
As in section \ref{sec:liqedge}, the covariance $\esp{\flucN{\e{e}}(\varphi_1)\flucN{\e{e}}S(\varphi_2)}$ is a convolution of distributions $\varphi_1^\epsi \ast F^\epsi$ applied to the test function $\varphi_2$, where

\begin{equation*}
\varphi_1^\epsi=\epsi^2\sum_{xy}\varphi_1(u_{\xx}) \delta(\cdot -
u^\epsi_{\xx}),\qquad
F^\epsi=\sum_{\xx}\mathrm{Cov}(\e{e},\e{e}_{\xx})\delta(\cdot-u^\epsi_{\xx}).
\end{equation*}

$\varphi_1^\epsi$ converges weakly to $\varphi_1$. We have now to prove that the distribution
$F^\epsi$ converges toward the distribution $F=\frac{\partial^2 \mathcal{F}}{\partial \log \K_{\e{e}}^2}\delta$. If it is the case, then the limit of the second moment will be
\begin{equation}
\langle \varphi_1 \ast F,\varphi_2\rangle=\frac{\partial^2 \mathcal{F}}{\partial \log \K_{\e{e}}^2}\langle \varphi_1\ast\delta,\varphi_2\rangle = \frac{\partial^2 \mathcal{F}}{{\partial \log \K_{\e{e}}}^2}\int_{\RR^2}\varphi_1(u)\varphi_2(u)|\ud u|.
\end{equation}

If $\xx\neq(0,0)$, the covariance between edges $\e{e}$ and $\e{e}_{\xx}$ is given by

\begin{equation}
\mathrm{Cov}(\e{e},\e{e}_{\xx})=\esp{(\e{e}_{0,0}-\overline{\e{e}})(\e{e}_{\xx}-\overline{\e{e}})}=-\K_{\e{e}}^2 \Ki(\xx)\Ki(-\xx)
\end{equation}

Let $\psi$ be a smooth test function with compact support, and $N$ a large integer. We decompose the sum over $\xx$ in the expression of $\langle F^\epsi,\psi\rangle$ depending on whether the norm of $\xx$ is larger than $N$ or not.

\begin{equation*}
\langle F^\epsi,\psi\rangle = \sum_{|\xx|>N}
\psi(u^\epsi_{\xx})\mathrm{Cov}(\e{e},\e{e}_{\xx}) +\sum_{|\xx|\leq
N} \left(\psi(u_{\xx}^\epsi)-\psi(0)\right)
\mathrm{Cov}(\e{e},\e{e}_{\xx})
 + \psi(0)\sum_{|\xx|\leq N}
\mathrm{Cov}(\e{e},\e{e}_{\xx}).
\end{equation*}

In the first sum, there are at most $O(\epsi^{-2})$ terms since the support of $\psi$ is bounded. As $|\xx|>N$, each term in this sum is bounded by some constant times $e^{-2C_2 N}$. Therefore the whole sum is a $O(\epsi^{-2}e^{-2C_2 N})$.
In the second sum, since $|\xx|\leq N$, the distance between $u^\epsi_{\xx}$ and $0$ is less than $\epsi N$. As $\psi$ is smooth , $\left(\psi(u^\epsi_{\xx})-\psi(0)\right)$ is $O(\epsi N)$. Since there are $O(N^2)$ terms in the second sum, it is a $O(N^3 \epsi)$.
If we choose for instance $N$ of order $\epsi^{-1/4}$, these two sums converge to zero, when $\epsi$ goes to zero.

The third sum is absolutely convergent, the limit of $\langle F^\epsi,\psi\rangle$ is
\begin{equation*}
\psi(0)\sum_{\xx\in\mathbb{\ZZ^2}}\mathrm{Cov}(\e{e},\e{e}_{\xx}),
\end{equation*}
meaning that $F^\epsi$ converges weakly to $\sum_{\xx\in\mathbb{\ZZ^2}}\mathrm{Cov}(\e{e},\e{e}_{\xx})\delta$ and the limit of the second moment is proportional to the $\mathrm{L}^2$ scalar product. The coefficient of proportionality
\begin{equation*}
\displaystyle\sum_{\xx\in\mathbb{\ZZ^2}}\mathrm{Cov}(\e{e},\e{e}_{\xx})
\end{equation*}
 can be rewritten in terms of the polynomials
$P(z,w)$ and $Q_{\bv\wv}(z,w)=Q_\e{e}(z,w)$.

\begin{align*}
\sum_{\xx\in\ZZ^2}\mathrm{Cov}(\e{e},\e{e}_{\xx})&=\esp{(\e{e}-\overline{\e{e}})(\e{e}-\overline{\e{e}})}+\sum_{\xx\neq(0,0)}\esp{(\e{e}_{0,0}-\overline{\e{e}})(\e{e}_{\xx}-\overline{\e{e}})}  \\
&=\prob{\e{e}}(1-\prob{\e{e}})-\K_\e{e}^2 \sum_{\xx\neq(0,0)}\Ki(\xx)\Ki(-\xx) \\
&=\K_\e{e}\Ki(0)-\K_\e{e}^2\sum_{\xx}\Ki(\xx)\Ki(-\xx)\\
&=\iint_{\mathbbm{T}^2} \frac{\K_\e{e}Q_\e{e}(z,w)}{P(z,w)}-\frac{\K_\e{e}^2 Q_\e{e}(z,w)^2}{P(z,w)^2}\frac{\ud z}{2\pi iz}\frac{\ud w}{2\pi iw}. 
\end{align*}
Since $Q_{\e{e}}(z,w)=\frac{\partial}{\partial \K_{\e{e}}} P(z,w)$, we have finally
\begin{align*}
\sum_{\xx\in\ZZ^2}\mathrm{Cov}(\e{e},\e{e}_{\xx})
&=\K_{\e{e}}\frac{\partial}{\partial \K_{\e{e}}}\K_{\e{e}}\frac{\partial}{\partial \K_{\e{e}}} \iint_{\mathbbm{T}^2} \log(P(z,w)) \frac{\ud z}{2\pi iz}\frac{\ud w}{2\pi iw}= \frac{\partial^2 \mathcal{F}}{\partial \log \K_{\e{e}}^2}.
\end{align*}\qed
\end{proof}

\subsection{Higher moments}

\begin{proposition}
The Wick formula is verified in the limit.
\begin{equation*}
\lim_{\epsi \rightarrow 0} \esp{(\flucN{\e{e}}\psi)^n)}=\left\{\begin{array}{cl} 0 & \textrm{if $n$ is odd}, \\ (n-1)!!\esp{(\flucN{\e{e}})^2}^{n/2} & \textrm{if $n$ is even}.\end{array}\right.
\end{equation*}
\end{proposition}

\begin{proof}
As in the case of a generic liquid measure, the proof begins with the study of the restricted
 $n$th moment $\Xi^\epsi_n(\psi)$ defined by equation \eqref{eq:Xi_n}.
  Lemma \ref{lem:corr} yields an asymptotic expression with the same structure as in the generic liquid case \eqref{eq:dvpmmt}:

\begin{multline}
\Xi_n^\epsi(\psi)=
\sum_{\{\Gamma\}_{l=1}^p}\prod_{l=1}^p\sum_{\substack{\gamma\textrm{
cycle}\\
\mathrm{supp}(\gamma)=\Gamma_l}}\biggl(\sgn(\gamma)(\epsi\K_{\e{e}})^{|\gamma|}\times\biggr.\\
\biggl.\sum_{\substack{\xx_{j_1},\cdots,\xx_{j_{|\gamma|}}
\\distinct}}\prod_{k=1}^{|\gamma|}
\psi(u^\epsi_{\xx_{j_k}})\Ki({\xx_{\gamma(j_k)}}-{\xx_{j_k}})
+o(1)\biggr)
\end{multline}

The contributions of cycles of length greater than 3 vanish in the limit, as an application of the following lemma:

\begin{lemma}
\begin{equation*}
\forall m\geq 3,\quad
 \lim_{\epsi\rightarrow 0}\epsi^m\sum_{\xx_1,\cdots,\xx_m}\prod_{j=1}^{m}\psi(u^{\epsi}_{\xx_j})\Ki(\xx_{j+1}-\xx_{j})=0
 \end{equation*}
\end{lemma}

\begin{proof}
Define $\xx'=\xx_1$ and for $j\geq 2$ $\xx'_j=\xx_j-\xx_{j-1}$ . The sum can be rewritten using these new notations and invariance by translation of operator $\Ki$ as

$$\epsi^{m}\sum_{\xx'}\psi(u^\epsi_{\xx'})\sum_{\xx'_2\cdots\xx'_m}\left(\prod_{j=2}^{m}\psi(u^{\epsi}_{\xx'+\xx'_2+\cdots+\xx'_j})\Ki(\xx'_j)\right)\Ki(-\xx'-\xx'_2\cdots-\xx'_m)$$

As $\psi$ is a continuous function on a compact set, it is bounded. The sum on $\xx'$ has $O(\epsi^{-2})$ terms. $\Ki(\bv_{-\xx'-\xx'_2\cdots-\xx'_n},\wv)$ is bounded independently from the $\xx'_j$'s. As $\Ki$ decays exponentially, the sum on $\xx'_2,\dots,\xx'_n$ is bounded. The whole sum is thus a $O(\epsi^{m-2})$, which goes to zero when $\epsi$ goes to zero as soon as $m\geq 3$.\qed
\end{proof}

Thus $\Xi_n^\epsi(\psi)$ converges to $\Xi_n(\psi)=(n-1)!!\ \Xi_2(\psi)$. The end of the proof deals with collisions between edges, which is identical to what has been done for proposition \ref{prop:nth}, leading to the result.\qed
\end{proof}

\subsection{Patterns in gaseous phase}\label{subsec:gaspat}

Combining the notations and the techniques introduced in section \ref{sec:liqpat} to deal with correlations between patterns, and following the steps of the proof of theorem \ref{thm:gasedge}, one can prove the following

\begin{theorem}
Let $\pat$ be a pattern in a dimer model endowed with a gaseous Gibbs measure. The random field $\flucN{\pat}$ of density fluctuation of pattern $\pat$ converges weakly in distribution to a white noise.
\end{theorem}
The proof is omitted here.

\section{Correlations between density fields\label{sec:corrfields}}
In the previous sections, only fluctuations of the density field associated to one fixed pattern was considered. One can ask what happens for correlations between density fields associated to different patterns. To what extent a high density of some pattern in a given region of the plane has an influence on the density of another pattern in an other region ?

This question is answered by the following theorem generalizing the results of the previous sections.
\begin{theorem}\label{thm:general}
 Consider a dimer model with a generic liquid Gibbs measure $\mu$.
\begin{itemize}
\item Let $\pat_1$ and $\pat_2$ be two patterns. The bilinear form $\esp{\flucN{\pat_1}(\cdot)\flucN{\pat_2}(\cdot)}$ on $\mathcal{C}^{\infty}_c(\RR^2)$ converges when $\epsi$ goes to zero to a bilinear form $\esp{\lflucN{\pat_1}(\cdot)\lflucN{\pat_2}(\cdot)}$.
\item
If $\mu$ is a generic liquid Gibbs measure, there exists a constant $A_{\pat_1\pat_2}$ such that for every test functions $\varphi_1$ and $\varphi_2$,
\begin{multline}
\esp{\lflucN{\pat_1}(\varphi_1)\lflucN{\pat_2}(\varphi_2)}=\frac{1}{\pi}\iint_{\RR^2\times\RR^2}\partial_{\pat_1^*}\varphi_1(u)G(u,v)\partial_{\pat_2^*}\varphi_2(v)|\ud
u||\ud v|\\+ A_{\pat_1\pat_2}
\int_{\RR^2}\varphi_1(u)\varphi_2(u)|\ud u|
\end{multline}
\item If $\mu$ is gaseous, there exists a constant $A_{\pat_1\pat_2}$ such that for every test functions $\varphi_1$ and $\varphi_2$,
\begin{equation*}
\esp{\lflucN{\pat_1}(\varphi_1)\lflucN{\pat_2}(\varphi_2)}= A_{\pat_1\pat_2} \int_{\RR^2}\varphi_1(u)\varphi_2(u)|\ud u|
\end{equation*}
\item Let $\pat_1,\dots,\pat_n$ be patterns (not necessarily distinct). When $\epsi$ goes to zero,
the multilinear form $\esp{\flucN{\pat_1}(\cdot)\cdots\flucN{\pat_n}(\cdot)}$ converges.
 The limit $\esp{\lflucN{\pat_1}(\cdot)\cdots\lflucN{\pat_n}(\cdot)}$ is given by Wick formula:
  for every test functions $\varphi_1,\dots,\varphi_n$,
\begin{equation*}
\esp{\lflucN{\pat_1}(\varphi_1)\cdots\lflucN{\pat_n}(\varphi_n)}=\begin{cases}
 0 & \text{if $n$ is odd} \\
\displaystyle\sum_{\text{pairings}} \prod_{k=1}^{n/2} \esp{\lflucN{\pat_{i_k}}(\varphi_{i_k})\lflucN{\pat_{j_k}}(\varphi_{j_k})} & \text{if $n$ is even}
\end{cases}
\end{equation*}
\end{itemize}
\end{theorem}
The coefficients $A_{\pat_1\pat_2}$ are in general not known in a closed form. However some relations between them can be found. For example, let $\mathbf{v}$ be a vertex of the graph, and denote by $\e{e}_1,\dots,\e{e}_m$ the edges incident with $\mathbf{v}$. The complex numbers $\e{e}_1^*,\dots,\e{e}_m^*$ sum to zero since they represent the edges of the dual face $\mathbf{v}^*$. Even more, for every $\epsi>0$, the linear combination $\flucN{\e{e}_1}+\cdots+\flucN{\e{e}_m}$ is identically zero. Indeed, for every $\xx\in\ZZ^2$, there is exactly one edge incident with $\mathbf{v}_\xx$ in the random dimer configuration. Therefore the sum of indicator functions $(\e{e}_1)_\xx+\cdots+(\e{e}_m)_\xx$ is always equal to 1. This relation between the random fields $\flucN{\e{e}_j}$ at a microscopic level yields relations between the different coefficients $A_{\e{e}_i\e{e}_j}$. Precisely,
\begin{equation}
\sum_{i=1}^m\sum_{j=1}^m A_{\e{e}_i\e{e}_j}=0.
\end{equation}

\section{Examples}

The theorems in the previous sections state a convergence of density fluctuation in the scaling limit to a linear combination of a derivative of the massless free field and a white noise. However, they do not give an explicit form for the white noise amplitude. 
In this section we present some cases for which a closed expression for the white noise amplitude can be provided in terms of the weights on edges.
The first case is the dimer model on the graph $\ZZ^2$ with periodic weights $a,b,c,d$ around white vertices.
The second case is the dimer model on the square-octagon graph.

\subsection{Dimer densities on $\ZZ^2$}

\begin{figure}
\centering
\includegraphics[height=5cm]{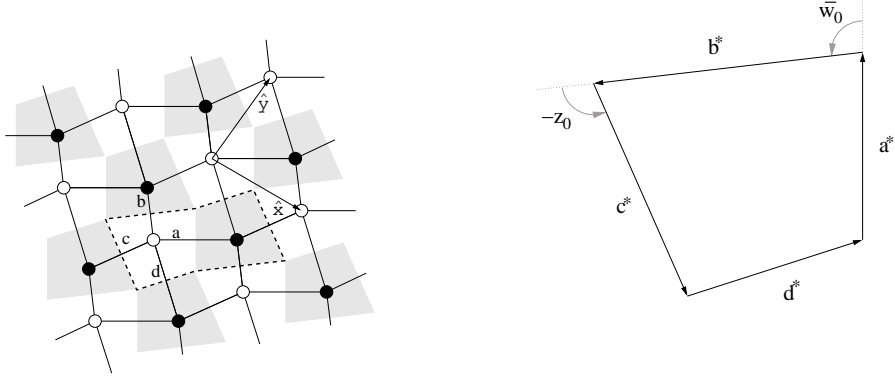}
\caption{\footnotesize On the left, a piece of the geometric realization (isoradial embedding in this case) of $\ZZ^2$ fixed by the weights $a,b,c,d$ assigned to edges. The region delimited by the tick dotted contour is a fundamental domain, and the coloring-preserving symmetries are generated by $\hat{\mathrm{x}}$ and $\hat{\mathrm{y}}$. On the right, a white face and the quantities (angles and side lengths) related to it.\label{fig:isoZ2}}
\end{figure}

The graph we consider here in the graph $\ZZ^2$ with a bipartite coloring of its vertices. Weights are assigned to edges according to their directions: $a,b,c,d$ counterclockwise around white vertices and clockwise around black vertices. If none of the weights is greater than the sum of the others, the corresponding dimer model is critical \cite{Ke:Crit}:  the graph can be embedded in the plane such that all the faces of the graph as well as those of the dual graph are inscribed in circles of a given radius (see figure~\ref{fig:isoZ2}). The Gibbs measure with no magnetic field on dimer configurations is liquid. The dual faces are similar to the cyclic quadrilateral with sides $a,b,c$ and $d$. The area of such a quadrilateral is
\begin{equation*}
\mathrm{Area}=\frac{1}{4}\sqrt{(-a+b+c+d)(a-b+c+d)(a+b-c+d)(a+b+c-d)}
\end{equation*}
and the radius $R$ of its circumscribed circle is defined by the relation
\begin{equation*}
R^2=\frac{(ab+cd)(ac+bd)(ad+cb)}{(-a+b+c+d)(a-b+c+d)(a+b-c+d)(a+b+c-d)}.
\end{equation*}

The fact we chose the fundamental domain to have area $1$ leads to the following expression for the complex numbers representing the dual edges in the embedding of $\ZZ^2$:
\begin{equation*}
\e{a}^*=\frac{ia}{\sqrt{\mathrm{2\ Area}}} \quad \e{b}^*=\frac{ib }{w_0\sqrt{2\ \mathrm{Area}}} \quad \e{c}^*=\frac{-ic z_0 }{w_0\sqrt{2\ \mathrm{Area}}} \quad \e{d}^*=\frac{id z_0}{\sqrt{2\ \mathrm{Area}}} \quad
\end{equation*}
where $(z_0,w_0)$ is the root of the the characteristic polynomial \begin{equation*}
P(z,w)=a+\frac{b}{w}-\frac{cz}{w}+dz
\end{equation*}
 on the unit torus, with the additional constraint that $\Im(z_0)>0$.

The following theorem is a particular case of theorem \ref{thm:general}.

\begin{theorem} Let $GFF$ a Gaussian free field in the plane and $W$ an independent white noise with unit variance.
The vector-valued random field
\begin{equation*}
\flucN{}=\begin{pmatrix}\flucN{a} \\ \flucN{b} \\ \flucN{c} \\
\flucN{d}\end{pmatrix}
\end{equation*}
 converges weakly in distribution to
the vector-valued Gaussian Field
\begin{equation*}
\begin{pmatrix} \lflucN{a} \\ \lflucN{b} \\ \lflucN{c} \\ \lflucN{d}\end{pmatrix}=
\frac{1}{\sqrt{\pi}}\begin{pmatrix} \partial_{\e{a}^*} \\ \partial_{\e{b}^*} \\ \partial_{\e{c}^*} \\ \partial_{\e{d}^*}\end{pmatrix}GFF
+\sqrt{\frac{abcd}{8\pi R^2 \ \mathrm{Area}}} \begin{pmatrix} +1 \\ -1 \\ +1 \\ -1\end{pmatrix}W.
\end{equation*}
\end{theorem}

\begin{proof}
The only point we have to explain is the computation of the white noise amplitude. This coefficient can be identified from the limit of the second moment of $\flucN{a}$. The proof of the convergence of the second moment, as before, goes through the proof of the convergence of the distribution $F^\varepsilon$ defined by \eqref{eq:Feps}. The sum over $\xx$ in the definition of $F^\epsi$ is decomposed into two parts depending on whether $u^\epsi_\xx$ belongs to some neighborhood of 0 or not. In the general case, we considered the neighborhood $\mathcal{B}=\left\{i(\alpha z_0 s +\beta z_0 t) \ ;\
(s,t)\in[-1,1]^2\right\}$. However, here we will take an infinite strip
\begin{equation*}
\mathcal{S}=\left\{i(\alpha z_0 s +\beta w_0 t) \ ;\
(s,t)\in\RR\times[-1,1]\right\}.
\end{equation*}
The condition on $\xx=(x,y)$ corresponding to $u^\epsi_\xx\in \mathcal{S}$ is $x\in\ZZ$ and $|y|\leq M$ where $M=\lfloor 1/\epsi\rfloor$. We have
\begin{equation*}
\langle F^\epsi,\psi\rangle = \sum_{\substack{ x\in \ZZ\\ |y|>M }}
\psi(u^\epsi_{\xx})\mathrm{Cov}(\e{e},\e{e}_{\xx})
+\sum_{\substack{ x\in \ZZ\\ |y|\leq M }}
\left(\psi(u^\epsi_{\xx})-\psi(0)\right)\mathrm{Cov}(\e{e},\e{e}_{\xx}) +\psi(0)\sum_{\substack{ x\in \ZZ\\ |y|\leq M }}
\mathrm{Cov}(\e{e},\e{e}_{\xx}).
\end{equation*}

The inverse Kasteleyn operator in this case is given by
\begin{equation*}
\Ki(x,y)=\Ki(\bv_\xx,\wv_0)=\iint_{\mathbbm{T}^2} \frac{z^{-y} w^{x}}{a+b/w -cz/w +dz} \frac{\ud z}{2i\pi z}\frac{\ud w}{2 i \pi w}.
\end{equation*}

The fact we have an infinite strip allows us to make use of one dimensional Fourier series in the $x$-direction to compute the third sum. Indeed, $\Ki(x,y)$ is the $x$th Fourier coefficient of the function $f_y(w)$ defined by
\begin{equation*}
f_y(w)=\int_{\cercle}\frac{z^{-y}}{a+b/w-cz/w+d z}\frac{\ud z}{2i\pi z}.
\end{equation*}
Hence, for a fixed $y$,
\begin{equation*}
\sum_{x\in\ZZ} \Ki(x,y)\Ki(-x,-y) = \int_{\cercle} f_y(w) f_{-y}(w) \frac{\ud w}{2i\pi w}.
\end{equation*}
For $y\neq0$, the functions $f_y$ and $f_{-y}$ have disjoint support, therefore the sum above is zero. When $y=0$, the sum is equal to
\begin{equation*}
\sum_{x\in\ZZ} \Ki(x,0)\Ki(-x,0) = \int_{\cercle} f_{0}^2(w) \frac{\ud w}{2i\pi w}=\int_{\left|\frac{a+bw}{cw-d}\right|>1} \frac{\ud w}{2i\pi(a+b/w)^2 w}.
\end{equation*}
The third sum equals
\begin{multline*}
\psi(0)\sum_{\substack{ x\in \ZZ\\ |y|\leq M }}
\mathrm{Cov}(\e{e},\e{e}_{\xx})=\psi(0)\left(a\Ki(0)-a^2\sum_{|y|\leq M}\sum_{x\in\ZZ}\Ki(x,y)\Ki(-x,-y)\right)\\
=\psi(0)\left(\int_{\left|\frac{a+bw}{cw-d}\right|>1} \frac{a\ud z}{(a+b/w)2i\pi w}-\int_{\left|\frac{a+bw}{cw-d}\right|>1} \frac{a^2\ud z}{(a+b/w)^2 2i\pi w}\right)
=\psi(0)\frac{-a b \Im(w_0)}{\pi|a+b/w_0|^2}.
\end{multline*}

The other term coming from the application of Green formula can also be computed and turns out to be equal to
\begin{equation*}
-\psi(0)\frac{1}{2\mathrm{Area}\ \pi}\left(\frac{-ab \Im(w_0)}{|a+\frac{b}{w_0}|}\right)^2.
\end{equation*}
Hence the variance of the white noise appearing in the limit of $\flucN{a}$ is given, after some calculations, by
\begin{equation*}
\frac{1}{\pi} \frac{-a b \Im(w_0)}{|a+b/w_0|^2}\left(1+\frac{ab \Im(w_0)}{2\mathrm{Area}}\right)= \frac{abcd}{8\pi R^2 \ \mathrm{Area}}
\end{equation*}

A similar computation for $\esp{\flucN{a}(\varphi_1)\flucN{b}(\varphi_2)}$ leads to the same expression, with a negative sign. As the expression of the coefficient is invariant under cyclic permutation of $a,b,c,d$, the amplitudes for the other pairs are easily deduced.\qed
\end{proof}


When $d=0$, the dimer model on $\ZZ^2$ is equivalent to that on the honeycomb lattice with periodic weights $a,b,c$. One can notice that in this case, the amplitude of the white noise vanishes. The interaction between dimers on the honeycomb lattice is purely electrostatic. We conjecture that it is true only for that particular model.

\subsection{Dimer densities on the square-octagon graph}

\begin{figure}
\begin{center}
\includegraphics[width=6cm]{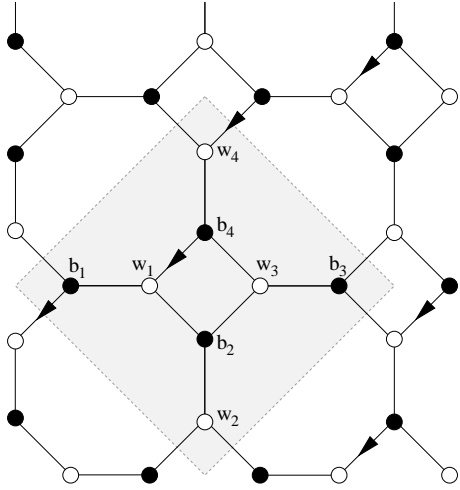}
\caption{\footnotesize A portion of the square-octagon graph. Edges are oriented from white end to black end, except those whose orientation is represented on the figure.}
\label{fig:squareoct}
\end{center}
\end{figure}

The square-octagon graph is a $\ZZ^2$-periodic graph whose fundamental domain is presented in figure \ref{fig:squareoct}. It contains four white and four black vertices. When every edge is assigned a weight equal to 1, the characteristic polynomial is given by

\begin{equation*}
P(z,w)=\det\begin{bmatrix} 1 & \frac{1}{w} & 0 & -\frac{1}{z} \\ 1 & 1 & 1 & 0 \\ 0 & z & 1& w \\ -1 & 0 & 1 & 1 \end{bmatrix} = 5- z-\frac{1}{z} -w-\frac{1}{w}
\end{equation*}

When the magnetic field is zero, or weak enough, the dimer model is in a gaseous phase. The fluctuations of the density field of an edge can be therefore computed by taking derivatives of the free energy of the system with respect to the weights, as explained in section \ref{sec:gaz}.


To compute for instance the amplitude of the white noise in the limit density of the edges $(\wv_1,\bv_1)$,  we assign  to these edges a weight $e^a$ and to the others a weight equal to 1, and compute the second derivative of the free energy $\mathcal{F}$ associated to this model with respect to $a$. For these new weights, the characteristic polynomial is now
\begin{equation*}
P_a(z,w)=\det\begin{bmatrix} e^a & \frac{1}{w} & 0 & -\frac{1}{z} \\ 1 & 1 & 1 & 0 \\ 0 & z & 1& w \\ -1 & 0 & 1 & 1 \end{bmatrix} = 4+e^a-e^a z-\frac{1}{z}-e^a w- \frac{1}{w} .
\end{equation*}
If $a$ is small enough, $P_a(z,w)$ has no zeros on the unit torus, and the dimer model is in still a gaseous phase in absence of magnetic field. In every point of this gaseous phase, the free energy is constant and given by
\begin{align*}
\mathcal{F_a}&=\iint_{\tore} \log P_a(z,w) \frac{\ud z}{2i\pi z}\frac{\ud w}{2i\pi w}\\
&=\log(4+e^a)+\iint_{[0,2\pi]^2}\log\left(1-\frac{1}{4+e^a}\left(e^{a+i\theta}+e^{-i\theta}+e^{a+i\phi}+e^{-i\phi}\right)\right)\frac{\ud \theta}{2\pi}\frac{\ud \phi}{2\pi}
\end{align*}
Performing the change of variables $\alpha=\frac{\theta+\phi}{2},\beta=\frac{\theta-\phi}{2}$ and moving the contour of integration over $\alpha$ from $[0,2\pi]$ to $[-ia/2,-ia/2+2\pi]$ using analyticity and periodicity in $\alpha$, we finally get an expression of $\mathcal{F}$ in terms of an absolutely convergent series:
\begin{align}
\mathcal{F}_a&=\log(4+e^a)+\iint_{[0,2\pi]^2} \log\left(1-\frac{4e^{a/2}}{4+e^a}\cos(\alpha)\cos(\beta)\right)\frac{\ud\alpha}{2\pi}\frac{\ud\beta}{2\pi}\\
&=\log(4+e^a)-\sum_{k=1}^{\infty}\frac{1}{k}\left(\frac{4e^{a/2}}{4+e^a}\right)^k\left(\int_0^{2\pi}\cos^k(\alpha)\frac{\ud\alpha}{2\pi}\right)\left(\int_0^{2\pi}\cos^k(\beta)\frac{\ud\beta}{2\pi}\right)\\
&=\log(4+e^a)-\sum_{k=1}^{\infty}\frac{1}{2k}\left(\frac{e^{a/2}}{4+e^a}\right)^{2k}\left(\frac{(2k)!}{(k!)^2}\right)^2
\end{align}
$\mathcal{F}_a$ can be expressed as the value of a certain generalized hypergeometric function. The Taylor expansion up to order $2$, involving the complete elliptic integrals $K$ and $E$
 \begin{equation*}
 \mathcal{F}_a=\mathcal{F}+\left(\frac{1}{2}-\frac{3 K\left(\frac{16}{25}\right)}{5\pi}\right)a+\left(\frac{K\left(\frac{16}{25}\right)-E\left(\frac{16}{25}\right)}{2\pi}\right)\frac{a^2}{2}+O(a^3)
 \end{equation*}
 gives information on the statistics of the copies of edge $(\wv_1,\bv_1)$. The constant coefficient is the free energy of the initial model, the coefficient of $a$ is the probability of $(\wv_1,\bv_1)$
 \begin{equation*}
 \prob{(\wv_1,\bv_1)}=\frac{1}{2}-\frac{3 K\left(\frac{16}{25}\right)}{5\pi},
 \end{equation*}
 and the coefficient of $a^2/2$ gives the amplitude of the white noise describing the scaling limit of the fluctuations of the number of edges $(\wv_1,\bv_1)$
 \begin{equation*}
 \lim_{\epsi\rightarrow 0}\esp{(\flucN{(\wv_1,\bv_1)}\varphi)^2} =
 \frac{K\left(\frac{16}{25}\right)-E\left(\frac{16}{25}\right)}{2\pi} \int_{\RR^2}\varphi(u)^2\ud u.
 \end{equation*}

Similarly, one can compute the probability of seeing an edge of a square, for example $(\wv_2,\bv_1)$, and the amplitude of the white noise:
\begin{align*}
 \prob{(\wv_2,\bv_1)}&=\frac{1}{4}+\frac{3 K\left(\frac{16}{25}\right)}{10\pi},\\
 \lim_{\epsi\rightarrow 0}\esp{(\flucN{(\wv_1,\bv_1)}\varphi)^2} &= \frac{2 K\left(\frac{16}{25}\right)}{5\pi} \int_{\RR^2}\varphi(u)^2\ud u.
 \end{align*}

The fact we see elliptic functions showing up is related to the fact that the spectral curve $\{P_a(z,w)=0\}$ in this case is genus-1 algebraic curve.

\begin{paragraph}{Acknowledgements.} We warmly thank Richard Kenyon for proposing to study pattern densities in dimer models. We are grateful to him for the many fruitful discussions.
This work has been done when the author was at Unversit\'e Paris-XI. The last part of writing this paper was done in a project at CWI,
financially supported by the Netherlands Organization for Scientific
Research (NWO).
\end{paragraph}

\bibliographystyle{hsiam}
\bibliography{cedric}

\end{document}